\documentclass[11pt]{amsart}
\usepackage[T1]{fontenc}
\usepackage{lmodern}
\usepackage[margin=1in]{geometry}
\usepackage{amssymb,mathtools}
\usepackage{microtype}
\usepackage[hidelinks]{hyperref}
\allowdisplaybreaks[2]
\numberwithin{equation}{section}

\newcommand{\HH}{\mathbb H}
\newcommand{\RR}{\mathbb R}
\newcommand{\ZZ}{\mathbb Z}
\newcommand{\TT}{\mathbb T}
\DeclareMathOperator{\tr}{tr}
\DeclareMathOperator{\osc}{osc}

\theoremstyle{plain}
\newtheorem{theorem}{Theorem}[section]
\newtheorem{lemma}[theorem]{Lemma}
\newtheorem{proposition}[theorem]{Proposition}
\newtheorem{corollary}[theorem]{Corollary}
\newtheorem{quotedtheorem}{Theorem}

\theoremstyle{definition}
\newtheorem{definition}[theorem]{Definition}
\newtheorem{example}[theorem]{Example}
\theoremstyle{remark}
\newtheorem{remark}[theorem]{Remark}

\title[Dirichlet problem on flat quaternionic quotients]
{The Dirichlet problem for quaternionic Monge--Amp\`ere equations
on flat quaternionic quotients}

\author[X. Ma]{Xi-Nan Ma}
\address{School of Mathematical Sciences,
University of Science and Technology of China, Hefei, China}
\email{xinan@ustc.edu.cn}

\author[J. Zhu]{Jingyong Zhu}
\address{School of Mathematical Sciences, Sichuan University,
Chengdu 610065, China}
\email{jzhu@scu.edu.cn}

\date{}

\begin{document}

\begin{abstract}

We prove smooth solvability of the Dirichlet problem for quaternionic Monge--Amp\`ere 
equations on smooth relatively compact domains in $\HH^m/\Gamma$, where $m\ge2$ and 
$\Gamma$ is a discrete additive
subgroup: for a smooth positive density, 
the existence of a smooth strictly plurisubharmonic subsolution  guarantees the existence 
and uniqueness of the  smooth admissible solution. The technique of the proof 
is a simple adaptation from the Ph.D. Thesis of Jingyong Zhu \cite{Zhu} in 2017, 
supervised by  Xinan  Ma, 
which is a joint work of Ma and Zhu. 
The claim \cite[Theorem A.1]{OpenAI} by OpenAI
on Sogge's local smoothing conjecture is a special case of our theorem.
\end{abstract}

\maketitle

\section{Introduction}

Monge--Amp\`ere equations play a central role in nonlinear elliptic
PDE and differential geometry. Their geometric significance is illustrated
by S.-T.~Yau's celebrated solution of the Calabi conjecture
\cite{Yau1978} through the complex Monge--Amp\`ere equation, and by Pogorelov \cite{Po78},
Cheng and Yau's fundamental work \cite{ChengYau1976} on the regularity
of solutions to the higher-dimensional Minkowski problem.
For the Dirichlet problem, the pioneering work of Caffarelli, Nirenberg
and Spruck \cite{CNS1} established the existence of a unique strictly
convex solution, smooth up to the boundary, on smooth bounded strictly convex domains
in $\RR^n$, for smooth positive densities
and smooth boundary values. B.~Guan and Spruck \cite{GuanSpruck1993}
extended this theory to general smooth bounded domains, without
assuming boundary convexity, under the hypothesis of a smooth strictly
locally convex strict subsolution attaining the prescribed boundary
values. The subsolution thus provides
a way to obtain boundary estimates beyond convex domains.

For complex Monge--Amp\`ere equations, Caffarelli, Kohn, Nirenberg
and Spruck \cite{CKNS} established smooth Dirichlet problem solvability on
smooth bounded strictly pseudoconvex domains in $\mathbb C^n$.
B.~Guan \cite{Guan, Guan2009} established
Dirichlet problem solvability on general smooth bounded domains in complex
Euclidean space under a subsolution hypothesis. 

In the quaternionic $\HH^m$ setting, Alesker \cite{Alesker1,Alesker2}
developed quaternionic pluripotential theory and proved smooth
Dirichlet solvability for the quaternionic Monge--Amp\`ere equation
on Euclidean balls. For the quaternionic Calabi-Yau problem and related
equations on compact manifolds, see
\cite{Alesker2013,AleskerShelukhin2013,AleskerVerbitsky2010,GentiliVezzoni2026}.

In his Thesis, for domains in $\HH^m$, Zhu \cite[Theorem 1.1]{Zhu} extended the
smooth solvability for the Dirichlet problem  from balls to arbitrary smooth bounded domains
under a subsolution hypothesis, thereby answering the question raised
by Alesker \cite{Alesker2} on smooth Dirichlet solvability on
quaternionic strictly pseudoconvex domains.

The recent OpenAI preprint \cite{OpenAI}  claims a proof of Sogge's local
smoothing conjecture for the Euclidean wave equation in every spatial
dimension $n\ge2$. A crucial analytic ingredient in its
argument is \cite[Theorem A.1]{OpenAI},  
the solvability of the Dirichlet problem for a quaternionic Monge--Amp\`ere equation, 
which provides the
potentials used in the local smoothing estimates. The claim \cite[Theorem A.1]{OpenAI} 
treats product domains
$(\RR^m/\Lambda)\times V$, where $\Lambda\subset\RR^m$ is a full
lattice and $V\Subset\RR^{2m}$ is a smooth bounded domain. It assumes
that $V$ admits a smooth defining function whose quaternionic Hessian
is positive definite up to the boundary.

By a simple adaptation from the  Thesis of Jingyong Zhu \cite{Zhu}
supervised by Xinan Ma, which is a joint work of Ma and Zhu,  
we prove the existence of the smooth solution for the Dirichlet  problem  on
general smooth relatively compact domains in flat quaternionic
quotients, assuming only the existence of a smooth strictly
plurisubharmonic subsolution.
We treat both prescribed positive densities $f(q)$ and positive
densities $f(q,u)$ that are nondecreasing in $u$.
The subsolution supplies the global test function for derivative
estimates, while boundary estimates are obtained in local lifts to
quaternionic space.  The claim \cite[Theorem A.1]{OpenAI} 
is a special case of our theorem \ref{thm:main}.
Our domain hypothesis is strictly more general: the periodic annuli
below admit suitable subsolutions but no strictly plurisubharmonic
defining function.

We first recall the periodic Dirichlet theorem of
\cite[Appendix A]{OpenAI} in our notation. Write
\[
q_a=x_a+iy_a+jz_a+kw_a,\qquad 1\le a\le m,
\]
and let $H$ denote the quaternionic Hessian and $\det$ its Moore determinant,
with the conventions in Section~\ref{sec:prelim}. Let $H_3$ be the restriction
of $H$ to functions independent of $w$. Explicitly,
\begin{equation}\label{eq:H3}
\begin{aligned}
H_3(u)_{ab}={}&u_{x_ax_b}+u_{y_ay_b}+u_{z_az_b}\\
&+i(u_{y_ax_b}-u_{x_ay_b})
 +j(u_{z_ax_b}-u_{x_az_b})
 +k(u_{z_ay_b}-u_{y_az_b}).
\end{aligned}
\end{equation}
In particular, $\tr H_3(u)=\Delta_{x,y,z}u$.

For a full lattice $\Lambda\subset\RR^m$, set
\[
Y=(\RR^m/\Lambda)\times\RR^{2m},
\]
where $x$ is periodic and $(y,z)\in\RR^{2m}$.

\begin{quotedtheorem}[\cite{OpenAI}]\label{thm:appendix}
Let $m\ge1$, let $V\Subset\RR^{2m}$ be a smooth domain, and set
$\Omega=(\RR^m/\Lambda)\times V$. Suppose that
$r\in C^\infty(\overline V)$ satisfies
\[
r<0\ \text{in }V,\qquad r=0,\ dr\ne0\ \text{on }\partial V,
\qquad H_3(r)>0\ \text{on }\overline\Omega.
\]
Let $B_0,f\in C^\infty(\overline\Omega)$ be real-valued, with
\[
H_3(B_0)\ge0,\qquad f>\det H_3(B_0)\quad\text{on }\overline\Omega.
\]
Then there is a unique $v\in C^\infty(\overline\Omega)$ satisfying
\[
\det H_3(B_0+v)=f,\qquad H_3(B_0+v)>0,\qquad
v|_{\partial\Omega}=0.
\]
Moreover, there are constants $b\ge b'>0$ such that
\begin{equation}\label{eq:two-sided}
br\le v\le b'r.
\end{equation}
For smooth finite-dimensional families of such data and domains,
the constants $b,b'$ can be chosen locally uniformly, and the solutions depend smoothly
on the parameters up to the boundary after smooth local identification of
the domains with a fixed domain.
\end{quotedtheorem}

For our nonlinear results, let $m\ge2$. Let $\Gamma$ be a discrete
additive subgroup of $\HH^m$ and set $X=\HH^m/\Gamma$, with its induced
flat metric and quaternionic structures. We consider a connected smooth
domain $M\Subset X$ with nonempty boundary. Functions on $M$ are pulled back by the quotient map to $\Gamma$-periodic functions on the preimage of $M$ in $\HH^m$.
A real function is called \emph{admissible} if $H(u)>0$. We first state
the result for a prescribed density.

\begin{theorem}\label{thm:main}
Let $f\in C^\infty(\overline M)$ be positive and let
$\phi\in C^\infty(\partial M)$. Suppose that there is a real function
$\underline u\in C^\infty(\overline M)$ satisfying
\begin{equation}\label{eq:subsolution}
H(\underline u)>0,\qquad
\det H(\underline u)\ge f\quad\text{in }\overline M,
\qquad \underline u=\phi\quad\text{on }\partial M.
\end{equation}
Then the Dirichlet problem
\begin{equation}\label{eq:main}
\begin{cases}
\det H(u)=f,\quad H(u)>0&\text{in }M,\\
u=\phi&\text{on }\partial M
\end{cases}
\end{equation}
has a unique admissible solution $u\in C^\infty(\overline M)$.
\end{theorem}

The determinant inequality in \eqref{eq:subsolution} need not be
strict. For a prescribed density $f(q)$, a smooth defining function $r$
with $H(r)>0$ on $\overline M$ supplies such a subsolution for every
smooth boundary value: if $\widetilde\phi$ is a smooth extension of
$\phi$, then $\underline u=\widetilde\phi+Ar$ satisfies
\eqref{eq:subsolution} for sufficiently large $A$. The subsolution
condition also applies to domains without such a defining function,
as the annuli below demonstrate.

Theorem~\ref{thm:solution-density} extends
Theorem~\ref{thm:main} to densities
depending on the solution. More precisely, if
$f\in C^\infty(\overline M\times\RR)$ satisfies $f>0$ and
$f_s=\partial f/\partial s\ge0$, the assumptions
\[
H(\underline u)>0,\qquad
\det H(\underline u)\ge f(q,\underline u),\qquad
\underline u|_{\partial M}=\phi
\]
guarantee a unique smooth admissible solution of
\[
\det H(u)=f(q,u)\quad\text{in }M,\qquad
u|_{\partial M}=\phi.
\]
 We prove the prescribed
density case first and give the changes in the global estimates for
$f(q,u)$ in Section~\ref{sec:solution-density}.

Theorem~\ref{thm:main} has the following consequence on the space $Y$
defined above.

\begin{corollary}\label{cor:three}
Let $\Omega\Subset Y$ be a smooth domain, $f\in C^\infty(\overline\Omega)$
be positive, and $\phi\in C^\infty(\partial\Omega)$.
If there is a smooth function $\underline u$ such that
\[
H_3(\underline u)>0,\qquad \det H_3(\underline u)\ge f,
\qquad \underline u|_{\partial\Omega}=\phi,
\]
then there is a unique smooth admissible solution of
\[
\det H_3(u)=f,\qquad H_3(u)>0,\qquad u|_{\partial\Omega}=\phi.
\]
\end{corollary}

Under the hypotheses of \cite[Theorem A.1]{OpenAI}, the function
$B_0+br$ is a subsolution for sufficiently large $b$. Consequently,
Corollary~\ref{cor:three} gives the existence assertion of
\cite[Theorem A.1]{OpenAI}  for $m\ge2$. The case $m=1$ is a linear
Dirichlet problem. The strict gap $f>\det H_3(B_0)$ provides a comparison
function $B_0+b'r$ for small $b'>0$, yielding the second inequality in
\eqref{eq:two-sided}. Section~\ref{sec:descent} recovers
\cite[Theorem A.1]{OpenAI}, including the locally uniform two-sided bounds
and smooth parameter dependence. The strict gap ensures
$v\le b'r<0$ in $\Omega$ but is not required for solvability in
Corollary~\ref{cor:three}.

The relaxation of the domain hypothesis is illustrated by the periodic annuli
\[
\Omega_0=\{(x,y,z)\in\TT^m\times\RR^{2m}:
                 a^2<|y|^2+|z|^2<b^2\},\qquad 0<a<b.
\]
Here $u_0=|y|^2+|z|^2$ satisfies $H_3(u_0)=4I_m$. Thus, for every
smooth density $0<f\le4^m$, Corollary~\ref{cor:three} gives a unique
smooth admissible solution with boundary values $a^2$ and $b^2$ on the
inner and outer components, respectively. These domains admit no
defining function with $H_3(r)>0$ on their closures: its quaternionic
tangential Hessian is negative on the inner boundary. Hence the
subsolution criterion covers Dirichlet data on domains excluded by \cite[Theorem A.1]{OpenAI}. Section~\ref{sec:examples} proves this
obstruction and gives a radial solution for constant densities.

We conclude the introduction with the main points of the estimates.
The Euclidean function $|q|^2$ need not descend to $X$, whereas the
subsolution is globally defined and satisfies
$H(\underline u)\ge\kappa I_m$ for some $\kappa>0$.
It supplies the global test function used to control first derivatives
and reduce the second derivative estimate to the boundary.
Constant real vector fields on $\HH^m$ descend to parallel fields on
$X$, so differentiation along these fields introduces no commutator
terms. Near the boundary, local inverses of the quotient map preserve
the quaternionic Hessian, and the local estimates of Ma-Zhu
\cite[Sections 3--4]{Zhu} apply. We include their detailed derivation
for prescribed densities. For $f(q,u)$, the nonpositive zeroth-order
term $-\partial_s\log f$ in the full linearization permits comparison
and the corresponding global derivative estimates. The boundary
estimates are again supplied by Ma-Zhu \cite{Zhu}.

More recently, using the technique of  the construction of suitable subsolutions, Ma and Zhang and their collaborators
\cite{MaZhangExterior,GMZRealPunctured,GMZComplexExterior,
	GMZComplexPunctured,GMYZ2026} studied Dirichlet problems for homogeneous
real and complex $k$-Hessian equations (the intermediate case between the Laplace and Monge--Amp\`ere equations) in exterior and punctured domains.
A common feature is the construction of suitable subsolutions and barriers
to obtain estimates for nondegenerate approximations that remain uniform
as the equations degenerate. These estimates lead to existence and
regularity results with prescribed behavior at infinity or near isolated
singularities.

The paper is organized as follows. Section~\ref{sec:prelim} fixes the quaternionic conventions and
collects the identities needed for linearization and comparison.
The a priori estimates occupy Sections~\ref{sec:estimates}
and~\ref{sec:normal}: the subsolution supplies global control
and reduces second derivative estimates to the boundary, where
the analysis is carried out in local quaternionic lifts.
These estimates yield smooth solvability and parameter dependence
in Section~\ref{sec:existence}, where we also prove
Corollary~\ref{cor:three}.
In Section~\ref{sec:descent}, we explain that \cite[Theorem A.1]{OpenAI} is a special case of Corollary~\ref{cor:three}.
The periodic annuli in Section~\ref{sec:examples} illustrate the
advantage of the subsolution hypothesis over a strictly
plurisubharmonic defining function.
Finally, Section~\ref{sec:solution-density} treats densities
$f(q,u)$, identifying the role of $f_u\ge0$ in comparison
and the global estimates.
\par

\begingroup
\subsection*{Acknowledgements} Xinan Ma was supported by the National Key Research and Development Program of China (Grant No. 2025YFA1017600) and NSFC No. 12141105. 

Statement of usage of artificial intelligence. The central mathematical ideas and overall
approach are due to the authors. OpenAI ChatGPT was used to assist with proofreading
and improving the clarity of the exposition. The authors take full responsibility for the
paper’s content and the accuracy of its references.
\leavevmode\par
\endgroup

\section{Quaternionic derivatives and linear algebra}\label{sec:prelim}

{
We fix the quaternionic notation and recall the basic algebraic and
analytic facts used in this paper. We also give proofs of the identities
and inequalities needed for the subsequent a priori estimates.
Unless otherwise indicated, the definitions and standard results follow
\cite[Section~1 and Subsection~2.1]{Alesker1}
and \cite[Section~2]{Zhu}.

\subsection{Quaternions and quaternionic space}\label{sec:quaternionic-space}

\begin{definition}[Quaternions]\label{def:quaternions}
The quaternions form the four-dimensional real algebra
\[
\HH=\{c=c_0+ic_1+jc_2+kc_3:c_0,c_1,c_2,c_3\in\RR\}.
\]
The coefficients $c_0,c_1,c_2,c_3$ are real and commute with $i,j,k$.
Addition is componentwise. Multiplication is real-bilinear and is
determined by
\begin{equation}\label{eq:quaternion-multiplication}
i^2=j^2=k^2=-1,\qquad
ij=-ji=k,\qquad jk=-kj=i,\qquad ki=-ik=j.
\end{equation}
\end{definition}

In particular, $ijk=-1$. For
$d=d_0+id_1+jd_2+kd_3$, expansion using \eqref{eq:quaternion-multiplication} gives
\[
\begin{aligned}
cd={}&(c_0d_0-c_1d_1-c_2d_2-c_3d_3)\\
&+i(c_0d_1+c_1d_0+c_2d_3-c_3d_2)\\
&+j(c_0d_2-c_1d_3+c_2d_0+c_3d_1)\\
&+k(c_0d_3+c_1d_2-c_2d_1+c_3d_0).
\end{aligned}
\]
Multiplication is associative and distributive, but is not commutative.
For example,
\[
(1+i)(1+j)=1+i+j+k,\qquad
(1+j)(1+i)=1+i+j-k.
\]
Thus one may change parentheses in a product, while keeping the
factors in their original order.

The conjugate, real part, and norm of $c$ are
\[
\bar c=c_0-ic_1-jc_2-kc_3,\qquad
\operatorname{Re}c=c_0=\frac{c+\bar c}{2},\qquad
|c|^2=c\bar c=\bar c c=\sum_{r=0}^3c_r^2.
\]
Conjugation reverses the order of a product:
\begin{equation}\label{eq:quaternion-conjugation}
\overline{cd}=\bar d\,\bar c,\qquad
\overline{c+d}=\bar c+\bar d,\qquad \overline{\bar c}=c.
\end{equation}
For instance, $\overline{ij}=-k=(-j)(-i)=\bar j\,\bar i$.
Every nonzero quaternion has a two-sided inverse,
\[
c^{-1}=\frac{\bar c}{|c|^2},\qquad cc^{-1}=c^{-1}c=1.
\]
The norm is multiplicative, as follows from
\[
|cd|^2=cd\bar d\bar c=|d|^2c\bar c=|c|^2|d|^2.
\]
Although $cd$ and $dc$ may differ, their real parts agree:
\begin{equation}\label{eq:quaternion-real-product}
\operatorname{Re}(cd)=c_0d_0-c_1d_1-c_2d_2-c_3d_3
=\operatorname{Re}(dc).
\end{equation}
By associativity, \eqref{eq:quaternion-real-product} gives cyclic invariance of the real part
of a product, such as $\operatorname{Re}(cde)=\operatorname{Re}(dec)$.

\begin{definition}[Quaternionic space and coordinates]\label{def:quaternionic-space}
The space $\HH^m$ consists of columns $q=(q_1,\ldots,q_m)^T$.
We use right scalar multiplication:
\[
\xi\alpha=(\xi_1\alpha,\ldots,\xi_m\alpha)^T,
\qquad \xi\in\HH^m,\quad \alpha\in\HH.
\]
Writing $q_a=x_a+iy_a+jz_a+kw_a$ identifies $\HH^m$ with $\RR^{4m}$ by
\begin{equation}\label{eq:real-coordinate-order}
q_{\RR}=(x_1,y_1,z_1,w_1,\ldots,x_m,y_m,z_m,w_m)^T,
\qquad |q|^2=\sum_{a=1}^m|q_a|^2.
\end{equation}
\end{definition}

For any quaternionic column $\xi$, the notation $\xi_{\RR}$ uses
the ordering in \eqref{eq:real-coordinate-order}. Open sets, smooth
functions, gradients, and real Hessians on $\HH^m$ are understood
through the identification $\HH^m\simeq\RR^{4m}$.

If $e_a$ is the column with $1$ in its $a$th entry and zero elsewhere,
then $\xi=\sum_a e_a\xi_a$.

\begin{definition}[Right-linear maps]\label{def:right-linear}
A map $T:\HH^m\to\HH^m$ is
right $\HH$-linear if
\[
T(\xi\alpha+\eta\beta)=(T\xi)\alpha+(T\eta)\beta.
\]
\end{definition}

Such a map is represented by a quaternionic matrix $A$ acting on the
left: the $b$th column of $A$ is $T(e_b)$, and
\begin{equation}\label{eq:matrix-action}
(A\xi)_a=\sum_{b=1}^m A_{ab}\xi_b,\qquad
A(\xi\alpha)=(A\xi)\alpha.
\end{equation}
The matrix coefficients precede the vector entries in each product.

\begin{definition}[Quaternionic affine lines]\label{def:quaternionic-line}
For $p\in\HH^m$ and $\xi\ne0$, a right quaternionic affine line is
\[
p+\xi\HH=\{p+\xi t:t\in\HH\}.
\]
\end{definition}

$p+\xi\HH$ is a real four-dimensional affine plane. Right and left scalar
multiplication need not give the same line. For example, with
$\xi=(1,j)^T\in\HH^2$,
\[
\xi i=(i,-k)^T,\qquad i\xi=(i,k)^T,
\]
and $i\xi\notin\xi\HH$: the first entry would force $t=i$, whereas
the second entry of $\xi t$ would then be $-k$.

\subsection{Hyperhermitian matrices and real representations}\label{sec:quaternionic-matrices}

For quaternionic matrices, multiplication is defined by
\[
(AB)_{ac}=\sum_b A_{ab}B_{bc}.
\]
Matrix multiplication is associative. The conjugate transpose is
\[
(A^*)_{ab}=\overline{A_{ba}},\qquad
\xi^*=(\bar\xi_1,\ldots,\bar\xi_m).
\]
The conjugation rule \eqref{eq:quaternion-conjugation} implies
\begin{equation}\label{eq:adjoint-identities}
(AB)^*=B^*A^*,\qquad (A^*)^*=A,
\qquad (A^{-1})^*=(A^*)^{-1}
\end{equation}
whenever $A$ is invertible. For example,
\[
((AB)^*)_{ca}
=\sum_b\overline{A_{ab}B_{bc}}
=\sum_b\overline{B_{bc}}\,\overline{A_{ab}}
=(B^*A^*)_{ca}.
\]
\begin{definition}[The standard quaternionic inner product]\label{def:quaternionic-inner-product}
For columns $\xi,\eta\in\HH^m$, whose entries $\xi_a,\eta_a$ are
quaternions, the standard quaternionic inner product is
\begin{equation}\label{eq:quaternionic-inner-product}
\langle\xi,\eta\rangle=\xi^*\eta=\sum_a\bar\xi_a\eta_a.
\end{equation}
\end{definition}

The inner product \eqref{eq:quaternionic-inner-product} satisfies
\[
\langle\xi,\xi\rangle=|\xi|^2,
\qquad \langle A\xi,\eta\rangle=\langle\xi,A^*\eta\rangle.
\]
Writing
\[
\xi_a=x_a+iy_a+jz_a+kw_a,\qquad
\eta_a=x'_a+iy'_a+jz'_a+kw'_a,
\]
we obtain
\begin{equation}\label{eq:real-inner-product}
\begin{aligned}
(\xi,\eta)_{\RR}
&:=\operatorname{Re}\langle\xi,\eta\rangle\\
&=\sum_{a=1}^m\operatorname{Re}\bigl[
(x_a-iy_a-jz_a-kw_a)(x'_a+iy'_a+jz'_a+kw'_a)\bigr]\\
&=\sum_{a=1}^m(x_ax'_a+y_ay'_a+z_az'_a+w_aw'_a)
=\xi_{\RR}^{T}\eta_{\RR}.
\end{aligned}
\end{equation}
Thus \eqref{eq:real-inner-product} identifies the real part of the
quaternionic inner product with the usual
Euclidean inner product on $\RR^{4m}$.
\begin{definition}[Unitary matrices]\label{def:unitary-matrix}
A matrix $U\in\HH^{m\times m}$ is called \emph{unitary} if $U^*U=I_m$.
The unitary matrices form the \emph{quaternionic unitary group}
$\operatorname{Sp}(m)$, also called the \emph{compact symplectic group}.
\end{definition}
Then $U^{-1}=U^*$ and $UU^*=I_m$. Equivalently, the columns
$v_1,\ldots,v_m$ of $U$ form an orthonormal quaternionic basis:
$v_a^*v_b=\delta_{ab}$, where $\delta_{ab}=1$ if $a=b$ and $0$ otherwise.
In particular,
\begin{equation}\label{eq:unitary-inner-product}
\langle U\xi,U\eta\rangle=\xi^*U^*U\eta=\langle\xi,\eta\rangle,
\qquad |U\xi|=|\xi|.
\end{equation}
By \eqref{eq:unitary-inner-product} and \eqref{eq:real-inner-product},
the elements of $\operatorname{Sp}(m)$ preserve the quaternionic inner
product and the Euclidean metric. Geometrically, quaternionic unitary
matrices describe changes of orthonormal quaternionic frames.

\begin{samepage}
\begin{definition}[Hyperhermitian semilinear forms]\label{def:hyperhermitian-form}
Let $V$ be a right quaternionic vector space. A map
$h:V\times V\to\HH$ is a \emph{hyperhermitian semilinear form}
if, for all vectors $\xi,\xi',\eta,\eta'\in V$ and all quaternionic
scalars $\alpha\in\HH$,
\begin{equation}\label{eq:hyperhermitian-form-properties}
\begin{aligned}
h(\xi+\xi',\eta)&=h(\xi,\eta)+h(\xi',\eta),\\
h(\xi,\eta+\eta')&=h(\xi,\eta)+h(\xi,\eta'),\\
h(\xi,\eta\alpha)&=h(\xi,\eta)\alpha,\\
h(\eta,\xi)&=\overline{h(\xi,\eta)}.
\end{aligned}
\end{equation}
Thus $h$ is additive in each variable, right-linear in the second
variable, and conjugate symmetric.
\end{definition}
\end{samepage}

For $\xi,\eta\in V$ and $\alpha,\beta\in\HH$,
\eqref{eq:hyperhermitian-form-properties} gives
\begin{equation}\label{eq:form-scalar-rule}
h(\xi\alpha,\eta\beta)
=\overline{h(\eta,\xi)\alpha}\,\beta
=\bar\alpha\,h(\xi,\eta)\beta.
\end{equation}
Thus $h$ is conjugate-linear in the first variable.
Conjugate symmetry also gives $h(\xi,\xi)\in\RR$.

\begin{definition}[Hyperhermitian matrices]\label{def:hyperhermitian-matrix}
A matrix $A\in\HH^{m\times m}$ is \emph{hyperhermitian} if
$A^*=A$, equivalently $A_{ba}=\overline{A_{ab}}$ for all $a,b$.
In particular, the diagonal entries $A_{aa}$ are real.
\end{definition}

\begin{proposition}[Matrix representation of a form]\label{prop:form-matrix}
Fix an ordered basis $e_1,\ldots,e_m$ of a right quaternionic vector
space $V$. Hyperhermitian semilinear forms on $V$ correspond bijectively
to hyperhermitian matrices $A\in\HH^{m\times m}$. In coordinate columns
$\xi,\eta\in\HH^m$, the correspondence is
\begin{equation}\label{eq:hyperhermitian-form}
h_A(\xi,\eta)=\xi^*A\eta
=\sum_{a,b}\bar\xi_a A_{ab}\eta_b,
\qquad A_{ab}=h(e_a,e_b).
\end{equation}
\end{proposition}

Taking $A=I_m$ recovers the standard inner product
\eqref{eq:quaternionic-inner-product}.

\begin{definition}[Positive forms and matrices]\label{def:positive-forms}
A hyperhermitian form $h$ is \emph{nonnegative definite} if
$h(\xi,\xi)\ge0$ for every $\xi\in V$, and \emph{positive definite}
if $h(\xi,\xi)>0$ for every nonzero $\xi\in V$.
For a hyperhermitian matrix $A$, the associated real quadratic form is
$Q_A(\xi):=h_A(\xi,\xi)=\xi^*A\xi$.
We write $A\ge0$ when $Q_A(\xi)\ge0$ for every $\xi\in\HH^m$,
and $A>0$ when $Q_A(\xi)>0$ for every $\xi\ne0$.
The order $A\ge B$ means $A-B\ge0$, and similarly for $A>B$.
\end{definition}

The reality of $Q_A$ can also be checked directly:
\[
\overline{\xi^*A\xi}=(\xi^*A\xi)^*=\xi^*A^*\xi=\xi^*A\xi.
\]

\begin{example}[Nonnegative matrices]\label{ex:nonnegative-matrices}
For columns $v,\xi\in\HH^m$ and a matrix $C\in\HH^{m\times m}$,
\[
\xi^*(vv^*)\xi=|v^*\xi|^2\ge0,
\qquad \xi^*(C^*C)\xi=|C\xi|^2\ge0.
\]
Thus $vv^*\ge0$ and $C^*C\ge0$.
If $C$ is invertible, then $C^*C>0$.
\end{example}

\begin{proposition}[Change of basis and positivity]\label{prop:form-change-basis}
Let $A$ be hyperhermitian and let $C\in\HH^{m\times m}$ be invertible.
Under the change of coordinate columns $\xi=C\widetilde\xi$, the
matrix of the form $h_A$ in the new basis is $C^*AC$.
The matrix $C^*AC$ is hyperhermitian, and
\begin{equation}\label{eq:positivity-congruence}
A>0\quad\Longleftrightarrow\quad C^*AC>0.
\end{equation}
Equivalence \eqref{eq:positivity-congruence} also holds with $\ge0$
in place of $>0$.
\end{proposition}

\begin{proof}
For coordinate columns $\xi=C\widetilde\xi$ and $\eta=C\widetilde\eta$
of two vectors, \eqref{eq:adjoint-identities} gives
\begin{equation}\label{eq:form-change-basis}
\xi^*A\eta=(C\widetilde\xi)^*A(C\widetilde\eta)
=\widetilde\xi^*(C^*AC)\widetilde\eta.
\end{equation}
Moreover, $(C^*AC)^*=C^*A^*C=C^*AC$.
To examine positivity, set $\widetilde\eta=\widetilde\xi$ in
\eqref{eq:form-change-basis}:
\begin{equation}\label{eq:quadratic-change-basis}
Q_A(\xi)=\xi^*A\xi
=\widetilde\xi^*(C^*AC)\widetilde\xi
=Q_{C^*AC}(\widetilde\xi).
\end{equation}
Since $C$ is invertible, $\xi=C\widetilde\xi$ ranges over all nonzero
columns in $\HH^m$ as $\widetilde\xi$ does. Formula
\eqref{eq:quadratic-change-basis} proves \eqref{eq:positivity-congruence}
and the nonnegative version.
\end{proof}

\begin{proposition}[Recovery from the quadratic form]\label{prop:quadratic-recovery}
If $A$ and $B$ are hyperhermitian and $Q_A(\xi)=Q_B(\xi)$ for all
$\xi\in\HH^m$, then $A=B$.
\end{proposition}

\begin{proof}
Let $e_1,\ldots,e_m$ be the standard coordinate columns in $\HH^m$.
The diagonal entries satisfy
\[
A_{aa}=Q_A(e_a)=Q_B(e_a)=B_{aa}.
\]
For $a<b$ and $t\in\HH$,
\begin{equation}\label{eq:quadratic-recovery}
Q_A(e_a+e_b t)=A_{aa}+|t|^2A_{bb}+2\operatorname{Re}(A_{ab}t).
\end{equation}
Write
\[
A_{ab}=c_0+ic_1+jc_2+kc_3,\qquad
B_{ab}=d_0+id_1+jd_2+kd_3,
\qquad c_r,d_r\in\RR.
\]
Taking $t=1,i,j,k$ in \eqref{eq:quadratic-recovery} gives
\begin{equation}\label{eq:quadratic-recovery-coefficients}
\begin{aligned}
c_0&=\operatorname{Re}A_{ab}
 =\tfrac12\bigl[Q_A(e_a+e_b)-Q_A(e_a)-Q_A(e_b)\bigr],\\
c_1&=-\operatorname{Re}(A_{ab}i)
 =\tfrac12\bigl[Q_A(e_a)+Q_A(e_b)-Q_A(e_a+e_b i)\bigr],\\
c_2&=-\operatorname{Re}(A_{ab}j)
 =\tfrac12\bigl[Q_A(e_a)+Q_A(e_b)-Q_A(e_a+e_b j)\bigr],\\
c_3&=-\operatorname{Re}(A_{ab}k)
 =\tfrac12\bigl[Q_A(e_a)+Q_A(e_b)-Q_A(e_a+e_b k)\bigr].
\end{aligned}
\end{equation}
Replacing $A$ by $B$ in \eqref{eq:quadratic-recovery-coefficients}
and using $Q_A=Q_B$ gives $c_r=d_r$ for $r=0,1,2,3$.
Hence, for $a<b$,
\[
A_{ab}=B_{ab},\qquad
A_{ba}=\overline{A_{ab}}=\overline{B_{ab}}=B_{ba}.
\]
Together with the diagonal equalities, we obtain $A=B$.
\end{proof}

For a square quaternionic matrix $A$, define
\[
\operatorname{Re}\tr A=\sum_a\operatorname{Re}A_{aa}.
\]
For compatible rectangular matrices $A,B$,
\eqref{eq:quaternion-real-product} gives
\begin{equation}\label{eq:real-trace-cyclic}
\operatorname{Re}\tr(AB)
=\sum_{a,b}\operatorname{Re}(A_{ab}B_{ba})
=\sum_{a,b}\operatorname{Re}(B_{ba}A_{ab})
=\operatorname{Re}\tr(BA).
\end{equation}
The real part in \eqref{eq:real-trace-cyclic} is essential.
For example, the one-by-one matrices $[i]$ and $[j]$ have products
$[k]$ and $[-k]$, whose traces differ.

\begin{definition}[The real representation of a quaternionic matrix]\label{def:real-representation}
In the real coordinates \eqref{eq:real-coordinate-order},
left multiplication by $c=c_0+ic_1+jc_2+kc_3$ has the matrix
\begin{equation}\label{eq:scalar-real-matrix}
R(c)=
\begin{pmatrix}
c_0&-c_1&-c_2&-c_3\\
c_1& c_0&-c_3& c_2\\
c_2& c_3& c_0&-c_1\\
c_3&-c_2& c_1& c_0
\end{pmatrix}
\end{equation}
in the ordered basis $(1,i,j,k)$. Thus the $4m\times4m$ real matrix
of $\xi\mapsto A\xi$ is the block matrix
\begin{equation}\label{eq:matrix-real-representation}
A_{\RR}=\bigl(R(A_{ab})\bigr)_{a,b=1}^m,
\qquad (A\xi)_{\RR}=A_{\RR}\xi_{\RR}.
\end{equation}
\end{definition}

\begin{proposition}[Compatibility with real linear algebra]\label{prop:real-representation}
The real representation preserves products and adjoints:
\begin{equation}\label{eq:real-matrix-identities}
(AB)_{\RR}=A_{\RR}B_{\RR},\qquad
(A^*)_{\RR}=A_{\RR}^T.
\end{equation}
If $A$ is hyperhermitian, then $A_{\RR}$ is real symmetric and
\begin{equation}\label{eq:real-quadratic-form}
\xi_{\RR}^T A_{\RR}\xi_{\RR}
=\operatorname{Re}(\xi^*A\xi)=\xi^*A\xi.
\end{equation}
Consequently, $A\ge0$ if and only if $A_{\RR}\ge0$.
If $U$ is quaternionic unitary, then $U_{\RR}$ is real orthogonal.
\end{proposition}

\begin{proof}
For every column $\xi\in\HH^m$, \eqref{eq:matrix-real-representation}
and associativity give
\[
(AB)_{\RR}\xi_{\RR}
=((AB)\xi)_{\RR}
=(A(B\xi))_{\RR}
=A_{\RR}(B\xi)_{\RR}
=A_{\RR}B_{\RR}\xi_{\RR}.
\]
Since $\xi_{\RR}$ ranges over $\RR^{4m}$, the product identity in
\eqref{eq:real-matrix-identities} follows.
The transpose identity follows from $R(\bar c)=R(c)^T$ in
\eqref{eq:scalar-real-matrix}. Apply \eqref{eq:real-inner-product}
to $\xi$ and $A\xi$ to obtain \eqref{eq:real-quadratic-form}.
For unitary $U$, \eqref{eq:real-matrix-identities} gives
\[
U_{\RR}^T U_{\RR}=(U^*U)_{\RR}=I_{4m}.
\]
\end{proof}

\begin{proposition}[Quaternionic spectral theorem]\label{prop:quaternionic-spectral}
Every hyperhermitian matrix $A$ admits a unitary diagonalization
\begin{equation}\label{eq:spectral-diagonalization}
U^*AU=\operatorname{diag}(\lambda_1,\ldots,\lambda_m),
\qquad U^*U=I_m,\qquad \lambda_a\in\RR.
\end{equation}
\end{proposition}

Writing the columns of $U$ as $v_1,\ldots,v_m$, the diagonalization means
\[
Av_a=v_a\lambda_a,\qquad
A=U\operatorname{diag}(\lambda_1,\ldots,\lambda_m)U^*.
\]
The eigenvalues are real: if $Av=v\lambda$ with $v\ne0$, then
$v^*Av=|v|^2\lambda$ is real. Moreover,
\[
\xi^*A\xi=\sum_a\lambda_a|(U^*\xi)_a|^2,
\]
so $A\ge0$ (respectively $A>0$) is equivalent to
$\lambda_a\ge0$ (respectively $\lambda_a>0$) for every $a$.
For $A>0$, the diagonalization \eqref{eq:spectral-diagonalization}
gives the inverse and the positive square root of $A$:
\begin{equation}\label{eq:spectral-inverse-square-root}
A^{-1}=U\operatorname{diag}(\lambda_a^{-1})U^*,\qquad
A^{1/2}=U\operatorname{diag}(\sqrt{\lambda_a})U^*.
\end{equation}

\subsection{Quaternionic derivatives and the Hessian}

\begin{definition}[Quaternionic differential operators]\label{def:quaternionic-derivatives}
For an open set $\Omega\subset\HH^m$ and $u\in C^1(\Omega,\HH)$,
define the quaternionic differential operators by
\begin{equation}\label{eq:quaternionic-derivatives}
\begin{aligned}
u_{\bar a}:=\frac{\partial u}{\partial\bar q_a}
&=u_{x_a}+iu_{y_a}+ju_{z_a}+ku_{w_a},\\
u_a:=\frac{\partial u}{\partial q_a}
&=\overline{\frac{\partial\bar u}{\partial\bar q_a}}
 =u_{x_a}-u_{y_a}i-u_{z_a}j-u_{w_a}k,
\end{aligned}
\end{equation}
where
\begin{equation}\label{eq:real-partial-derivatives}
u_{x_a}:=\frac{\partial u}{\partial x_a},\qquad
u_{y_a}:=\frac{\partial u}{\partial y_a},\qquad
u_{z_a}:=\frac{\partial u}{\partial z_a},\qquad
u_{w_a}:=\frac{\partial u}{\partial w_a}.
\end{equation}
\end{definition}

The quaternionic units $i,j,k$ act by left multiplication
in the first operator and by right multiplication in the second.
The operators in \eqref{eq:quaternionic-derivatives} are first-order
real differential operators with quaternionic coefficients. In one
variable, \eqref{eq:quaternionic-derivatives} gives
\begin{equation}\label{eq:coordinate-derivative-values}
\frac{\partial q}{\partial q}=4,\qquad
\frac{\partial q}{\partial\bar q}=-2,\qquad
\frac{\partial\bar q}{\partial\bar q}=4,\qquad
\frac{\partial\bar q}{\partial q}=-2.
\end{equation}
For example, $\partial q/\partial\bar q=1+i^2+j^2+k^2=-2$.
No factor $1/2$ is included in \eqref{eq:quaternionic-derivatives}.

\begin{definition}[Quaternionic Hessian]\label{def:quaternionic-hessian}
For a real function $u\in C^2(\Omega)$, the quaternionic Hessian is
\begin{equation}\label{eq:hessian}
\begin{aligned}
H(u)_{ab}=u_{\bar a b}
&:=\frac{\partial^2u}{\partial\bar q_a\,\partial q_b}
 =\frac{\partial}{\partial\bar q_a}
       \left(\frac{\partial u}{\partial q_b}\right)\\
&=u_{x_ax_b}+u_{y_ay_b}+u_{z_az_b}+u_{w_aw_b}\\
&\quad+i(u_{y_ax_b}-u_{x_ay_b}-u_{z_aw_b}+u_{w_az_b})\\
&\quad+j(u_{z_ax_b}-u_{x_az_b}+u_{y_aw_b}-u_{w_ay_b})\\
&\quad+k(u_{w_ax_b}-u_{x_aw_b}-u_{y_az_b}+u_{z_ay_b}).
\end{aligned}
\end{equation}
\end{definition}

\begin{remark}[Commutation of quaternionic derivatives]\label{rem:quaternionic-commutation}
For real coordinates $s,t$ and constant quaternions $c,d$,
commutation of real partial derivatives gives
\[
\begin{aligned}
c\,\frac{\partial}{\partial s}
  \left(\frac{\partial u}{\partial t}d\right)
&=c\,\frac{\partial^2u}{\partial s\,\partial t}d
 =c\,\frac{\partial^2u}{\partial t\,\partial s}d\\
&=\frac{\partial}{\partial t}
  \left(c\,\frac{\partial u}{\partial s}\right)d.
\end{aligned}
\]
Applying this identity to each term in
\eqref{eq:quaternionic-derivatives} yields
\begin{equation}\label{eq:quaternionic-commutation}
\frac{\partial}{\partial\bar q_a}
  \left(\frac{\partial u}{\partial q_b}\right)
=\frac{\partial}{\partial q_b}
  \left(\frac{\partial u}{\partial\bar q_a}\right).
\end{equation}
\end{remark}

Only quaternionic coordinate derivatives are abbreviated by $a$ or
$\bar a$. Real coordinate names are retained in subscripts, as in
$u_{x_a}$, $u_{x_ay_b}$, and $u_{x_a\bar b}=(u_{x_a})_{\bar b}$.
For real $u$, we have $u_a=\overline{u_{\bar a}}$.

\begin{proposition}[Hermitian symmetry and trace]\label{prop:hessian-trace}
For every real function $u\in C^2(\Omega)$, the matrix $H(u)$ is
hyperhermitian and
\begin{equation}\label{eq:trace}
H(u)_{aa}=u_{x_ax_a}+u_{y_ay_a}+u_{z_az_a}+u_{w_aw_a},
\qquad \tr H(u)=\Delta u.
\end{equation}
\end{proposition}

\begin{proof}
Since $u$ is real, \eqref{eq:quaternionic-derivatives} and
\eqref{eq:quaternionic-commutation} give
\begin{equation}\label{eq:hessian-conjugate-symmetry}
\overline{H(u)_{ab}}
=\overline{\frac{\partial}{\partial\bar q_a}
  \left(\frac{\partial u}{\partial q_b}\right)}
=\frac{\partial}{\partial q_a}
  \left(\overline{\frac{\partial u}{\partial q_b}}\right)
 =\frac{\partial}{\partial q_a}
  \left(\frac{\partial u}{\partial\bar q_b}\right)
=\frac{\partial}{\partial\bar q_b}
  \left(\frac{\partial u}{\partial q_a}\right)
 =H(u)_{ba}.
\end{equation}
Setting $b=a$ in \eqref{eq:hessian} gives the diagonal formula in
\eqref{eq:trace}, and summing over $a$ gives $\tr H(u)=\Delta u$.
\end{proof}

\begin{example}[Quadratic and affine functions]\label{ex:quadratic-functions}
For the real function $u(q)=|q|^2$,
\begin{equation}\label{eq:squared-norm-hessian}
\begin{aligned}
u_{\bar a}&=2q_a,\qquad u_b=2\bar q_b,\\
u_{\bar a b}&=\frac{\partial u_b}{\partial\bar q_a}
=2\frac{\partial\bar q_b}{\partial\bar q_a}
=8\delta_{ab},\qquad H(|q|^2)=8I_m.
\end{aligned}
\end{equation}
Here $\partial\bar q_b/\partial\bar q_a=4\delta_{ab}$ by
\eqref{eq:coordinate-derivative-values}.
More generally,
\[
H\left(\sum_{a=1}^m\mu_a|q_a|^2\right)
=8\operatorname{diag}(\mu_1,\ldots,\mu_m),\qquad \mu_a\in\RR,
\]
and every real affine function has zero quaternionic Hessian.
\end{example}

\begin{proposition}[Product and chain rules]\label{prop:hessian-product-chain}
For smooth real-valued functions $g,h$ and a smooth function
$\chi:\RR\to\RR$, the product and chain rules give
\begin{equation}\label{eq:product}
(gh)_{\bar a b}=g h_{\bar a b}+h g_{\bar a b}
                  +g_{\bar a}h_b+h_{\bar a}g_b.
\end{equation}
\begin{equation}\label{eq:chain-rule}
(\chi(g))_{\bar a b}=\chi'(g)g_{\bar a b}
                         +\chi''(g)g_{\bar a}g_b.
\end{equation}
\end{proposition}

\begin{proof}
The first derivatives are
\[
(gh)_b=g h_b+h g_b,\qquad
(\chi(g))_b=\chi'(g)g_b.
\]
Applying $\partial/\partial\bar q_a$ gives
\eqref{eq:product}--\eqref{eq:chain-rule}. The real-valued factors
$g,h$, and $\chi'(g)$ commute with the quaternionic coefficients.
The factors in the cross terms retain their displayed order.
\end{proof}

\subsection{Plurisubharmonic functions and domains}

\begin{definition}[Subharmonic functions]\label{def:subharmonic}
A real function $v$ on an open subset of $\RR^4$ is
subharmonic if $v$ is upper semicontinuous and satisfies the spherical
mean inequality
\[
v(t_0)\le\frac{1}{|\partial B_r|}
\int_{\partial B_r(t_0)}v\,dS
\]
for every closed ball contained in the domain of $v$.
For $v\in C^2$,
subharmonicity is equivalent to $\Delta v\ge0$.
\end{definition}

\begin{definition}[Quaternionic plurisubharmonic functions]\label{def:psh}
Let $\Omega\subset\HH^m$ be a domain. A real-valued function
$u:\Omega\to\RR$ is called \emph{quaternionic plurisubharmonic}
if $u$ is upper semicontinuous and the restriction of $u$ to every right
quaternionic line is subharmonic. A $C^2$ function $u$ is called
\emph{strictly plurisubharmonic} if the restriction of $u$ to every right
quaternionic line is strictly subharmonic, that is, has strictly
positive Laplacian.
\end{definition}

\begin{proposition}[Hessian criterion]\label{prop:psh-hessian}
Let $u\in C^2(\Omega)$ be real-valued. Then
\begin{equation}\label{eq:psh-hessian-criterion}
\begin{aligned}
u\text{ is quaternionic plurisubharmonic}
&\iff \xi^*H(u)(p)\xi\ge0
      \quad\text{for all }p\in\Omega,\ |\xi|=1\\
&\iff H(u)\ge0\quad\text{on }\Omega.
\end{aligned}
\end{equation}
Similarly, strict plurisubharmonicity is equivalent to $H(u)>0$.
\end{proposition}

\begin{proof}
For $p\in\Omega$ and $0\ne\xi\in\HH^m$, write
\begin{equation}\label{eq:quaternionic-line-parametrization}
u_{p,\xi}(t)=u(p+\xi t),\qquad
t=t_0+it_1+jt_2+kt_3,
\end{equation}
for $t$ near zero.
Let $A\in\HH^{m\times m}$ have first column $\xi$ and all other
columns zero. For a quaternionic column $q=(q_1,\ldots,q_m)^T$
near the origin, define
\begin{equation}\label{eq:line-extension}
\begin{aligned}
A&=(\xi\;0\;\cdots\;0),\qquad Aq=\xi q_1,\\
v(q)&=u(p+Aq)=u_{p,\xi}(q_1).
\end{aligned}
\end{equation}
The Hessian transformation formula in
\cite[Subsection~2.1, Claim~(i) preceding Proposition~2.1.6]{Alesker1}
applies to arbitrary quaternionic linear maps. Since translation
does not change the constant-coefficient differential operators in
\eqref{eq:quaternionic-derivatives}, the formula gives
\begin{equation}\label{eq:line-extension-hessian}
H(v)(0)=A^*H(u)(p)A.
\end{equation}
Because $A_{a1}=\xi_a$, the $(1,1)$ entry of
\eqref{eq:line-extension-hessian} is
\begin{equation}\label{eq:line-extension-entry}
H(v)(0)_{11}
=\sum_{a,b=1}^m\overline{A_{a1}}\,u_{\bar a b}(p)A_{b1}
=\sum_{a,b=1}^m\bar\xi_a\,u_{\bar a b}(p)\xi_b
=\xi^*H(u)(p)\xi.
\end{equation}
On the other hand, $v(q)=u_{p,\xi}(q_1)$ depends only on $q_1$.
Identifying $q_1$ with the quaternionic variable $t$ in
\eqref{eq:quaternionic-line-parametrization}, the definition
\eqref{eq:hessian} and the one-variable case of \eqref{eq:trace} yield
\begin{equation}\label{eq:line-restriction-hessian}
\begin{aligned}
\xi^*H(u)(p)\xi
&=H(v)(0)_{11}
 =\left.\frac{\partial}{\partial\bar q_1}
  \left(\frac{\partial v}{\partial q_1}\right)\right|_{q=0}\\
&=\left.\frac{\partial}{\partial\bar t}
  \left(\frac{\partial u_{p,\xi}}{\partial t}\right)\right|_{t=0}\\
&=\sum_{r=0}^3\frac{\partial^2u_{p,\xi}}{\partial t_r^2}(0)
 =\Delta_t u_{p,\xi}(0),\qquad
\Delta_t=\sum_{r=0}^3\frac{\partial^2}{\partial t_r^2}.
\end{aligned}
\end{equation}

When $|\xi|=1$, the map $t\mapsto p+\xi t$ is an isometry onto
$p+\xi\HH$, so $\Delta_t u_{p,\xi}(0)$ is the Laplacian of the
restriction of $u$ to that line at $p$.
Definition~\ref{def:psh} and
\eqref{eq:line-restriction-hessian} therefore give
\eqref{eq:psh-hessian-criterion}. The strict criterion follows
by replacing nonnegative Laplacians with positive Laplacians.
\end{proof}

\begin{example}[Comparison with convexity and subharmonicity]\label{ex:psh-comparison}
For a real $u\in C^2(\Omega)$, fix $p\in\Omega$ and
$0\ne\xi\in\HH^m$.
Let $D^2u(p)$ denote the real Hessian, with quaternionic directions
identified with real vectors as in \eqref{eq:real-coordinate-order}.
Since
$p+\xi t=p+t_0\xi+t_1\xi i+t_2\xi j+t_3\xi k$,
the real chain rule and \eqref{eq:line-restriction-hessian} give
\begin{equation}\label{eq:compression}
\begin{aligned}
\xi^*H(u)(p)\xi
&=\Delta_t u_{p,\xi}(0)\\
&=\left.\frac{d^2}{ds^2}\right|_{s=0}
  \bigl(u(p+s\xi)+u(p+s\xi i)+u(p+s\xi j)+u(p+s\xi k)\bigr)\\
&=D^2u(p)[\xi,\xi]+D^2u(p)[\xi i,\xi i]\\
&\quad+D^2u(p)[\xi j,\xi j]+D^2u(p)[\xi k,\xi k].
\end{aligned}
\end{equation}

When $m=1$, a right quaternionic line is all of $\HH$, so the
definition reduces to ordinary subharmonicity on $\RR^4$.
For $C^2$ functions in arbitrary dimension, \eqref{eq:compression}
and \eqref{eq:trace} give the implications
\begin{equation}\label{eq:convex-psh-subharmonic}
D^2u\ge0\quad\Longrightarrow\quad H(u)\ge0
\quad\Longrightarrow\quad\Delta u\ge0.
\end{equation}
By \eqref{eq:psh-hessian-criterion}, the three inequalities in \eqref{eq:convex-psh-subharmonic}
correspond to convexity, quaternionic plurisubharmonicity, and
subharmonicity, respectively, on a convex domain.
The converse implications fail in general. For example,
\[
u(q)=-x_1^2+y_1^2+z_1^2+w_1^2+\sum_{a=2}^m|q_a|^2
\]
has $H(u)=\operatorname{diag}(4,8,\ldots,8)>0$, although
$u_{x_1x_1}=-2$. Thus strict quaternionic plurisubharmonicity does not
imply real convexity. On $\HH^2$, the function
\[
v(q)=-|q_1|^2+2|q_2|^2
\]
satisfies $\Delta v=8>0$, but $H(v)=\operatorname{diag}(-8,16)$.
The restriction of $v$ to the line $q_2=0$ is not subharmonic.
\end{example}

\begin{definition}[Strictly pseudoconvex domains]\label{def:strictly-pseudoconvex}
A domain $\Omega\subset\HH^m$ with smooth boundary is called
\emph{strictly pseudoconvex} if, for every $p\in\partial\Omega$, there
exist a neighborhood $U$ of $p$ and a smooth strictly plurisubharmonic
function $\rho$ on $U$ such that
\[
\Omega\cap U=\{q\in U:\rho(q)<0\},\qquad
\rho(p)=0,\qquad d\rho(p)\ne0.
\]
\end{definition}

\begin{example}[Euclidean balls]\label{ex:strictly-pseudoconvex-ball}
The ball $B_R(p)$ has defining function
$\rho(q)=|q-p|^2-R^2$, for which
\[
H(\rho)=8I_m,\qquad |\nabla\rho|=2R\quad\text{on }\partial B_R(p).
\]
\end{example}

Definitions~\ref{def:psh} and~\ref{def:strictly-pseudoconvex} are
local and do not require boundedness.

On the quotient $X=\HH^m/\Gamma$, both definitions are applied in
local quaternionic lifts. Two such lifts differ by a translation,
which preserves $H(u)$ because the differential operators in
\eqref{eq:quaternionic-derivatives} have constant coefficients. Hence the definitions are independent of the
lift, even when the lifted domain is unbounded. The coordinate
transformation formulas are established in
Subsection~\ref{sec:quaternionic-coordinates}.

\subsection{The Moore determinant}
The independent real coordinates of a hyperhermitian matrix $A$
are the diagonal entries $A_{aa}$ and the four real components
of each $A_{ab}$ with $a<b$.

\begin{proposition}[Polynomial characterization]\label{prop:moore-polynomial}
There is a unique real polynomial $P$ in the independent real
coordinates of hyperhermitian $m\times m$ matrices such that
\begin{equation}\label{eq:moore-polynomial-characterization}
P(A)^4=\det_{\RR}A_{\RR}\quad\text{for every hyperhermitian }A,
\qquad P(I_m)=1.
\end{equation}
The polynomial $P$ is homogeneous of degree $m$ and has integer
coefficients.
\end{proposition}

\begin{definition}[Moore determinant]\label{def:moore-determinant}
For a hyperhermitian matrix $A$, the \emph{Moore determinant} is
\begin{equation}\label{eq:moore-definition}
\det A:=P(A),
\end{equation}
where $P$ is the polynomial in
Proposition~\ref{prop:moore-polynomial}.
\end{definition}

By \eqref{eq:moore-polynomial-characterization},
\begin{equation}\label{eq:moore-real-determinant}
(\det A)^4=\det_{\RR}A_{\RR}.
\end{equation}
Every monomial of $P$ has total degree $m$. In particular,
$\det(tA)=t^m\det A$ for $t\in\RR$.
The polynomial characterization fixes the sign of $\det A$, whereas
\eqref{eq:moore-real-determinant} alone determines only $|\det A|$.

\begin{lemma}[Moore determinant identities]\label{lem:moore-identities}
The Moore determinant agrees with the usual determinant on real
symmetric and complex Hermitian matrices.
For a hyperhermitian $A$ and any $C\in\HH^{m\times m}$, the
\emph{congruence identity} is
\begin{equation}\label{eq:moore-congruence}
\det(C^*AC)=\det A\,\det(C^*C).
\end{equation}
In particular, if $U\in\HH^{m\times m}$ satisfies $U^*U=I_m$,
then taking $C=U$ in \eqref{eq:moore-congruence} and using
$\det I_m=1$ gives
\begin{equation}\label{eq:moore-unitary}
\det(U^*AU)=\det A\,\det(U^*U)=\det A.
\end{equation}
Consequently, if $\lambda_1,\ldots,\lambda_m$ are the real diagonal
entries in \eqref{eq:spectral-diagonalization}, then
\begin{equation}\label{eq:moore-spectral}
\det A=\prod_{a=1}^m\lambda_a.
\end{equation}
\end{lemma}

\begin{proof}
For real polynomials $f,g$ in the same variables with
$f(x_0)=g(x_0)=1$,
\begin{equation}\label{eq:polynomial-fourth-power-uniqueness}
f^4\equiv g^4
\quad\Longrightarrow\quad
(f-g)(f+g)(f^2+g^2)\equiv0
\quad\Longrightarrow\quad f\equiv g.
\end{equation}
For complex Hermitian $A$, the ordinary complex spectral theorem
gives a complex unitary matrix $V$ such that
\[
V^*AV=\operatorname{diag}(\mu_1,\ldots,\mu_m),
\qquad V^*V=I_m,\qquad \mu_a\in\RR.
\]
By \eqref{eq:matrix-real-representation} and
\eqref{eq:real-matrix-identities},
\begin{equation}\label{eq:complex-spectral-real-representation}
V_{\RR}^{T}A_{\RR}V_{\RR}
=\operatorname{diag}(\mu_1I_4,\ldots,\mu_mI_4),
\qquad V_{\RR}^{T}V_{\RR}=I_{4m}.
\end{equation}
Taking real determinants and using \eqref{eq:moore-real-determinant} gives
\begin{equation}\label{eq:moore-complex-fourth-power}
(\det A)^4=\det_{\RR}A_{\RR}
=\prod_{a=1}^m\mu_a^4=(\det_{\mathbb C}A)^4.
\end{equation}
Both determinants are real polynomials in the independent real
coordinates of $A$ and equal $1$ at $A=I_m$, so
\eqref{eq:polynomial-fourth-power-uniqueness} gives their equality.
Real symmetric matrices are included.

By \eqref{eq:moore-real-determinant} and
\eqref{eq:real-matrix-identities},
\begin{equation}\label{eq:moore-congruence-fourth-power}
\begin{aligned}
\bigl(\det(C^*AC)\bigr)^4
&=\det_{\RR}(C_{\RR}^{T}A_{\RR}C_{\RR})\\
&=\det_{\RR}A_{\RR}\,\det_{\RR}(C_{\RR}^{T}C_{\RR})\\
&=\bigl(\det A\,\det(C^*C)\bigr)^4.
\end{aligned}
\end{equation}
The two sides of \eqref{eq:moore-congruence} are real polynomials in
the independent real coordinates of $A$ and the real coordinates
of $C$. At $A=C=I_m$,
$\det(C^*AC)=\det A\,\det(C^*C)=1$. Hence
\eqref{eq:polynomial-fourth-power-uniqueness} proves
\eqref{eq:moore-congruence}. Taking $C=U$ with $U^*U=I_m$ gives
\eqref{eq:moore-unitary}.
Finally, applying \eqref{eq:moore-unitary} to
\eqref{eq:spectral-diagonalization} gives
\[
\det A=\det(U^*AU)
=\det_{\RR}\operatorname{diag}(\lambda_1,\ldots,\lambda_m)
=\prod_{a=1}^m\lambda_a.\qedhere
\]
\end{proof}

Throughout the paper, $\det$ on a hyperhermitian matrix denotes
the Moore determinant.

\begin{example}[$2\times2$ matrices]\label{ex:moore-two-by-two}
For a $2\times2$ hyperhermitian matrix, with $r,s\in\RR$ and
$c=c_0+ic_1+jc_2+kc_3$, one has
\begin{equation}\label{eq:moore-two-by-two}
\det\begin{pmatrix}r&c\\\bar c&s\end{pmatrix}
=rs-|c|^2=rs-c_0^2-c_1^2-c_2^2-c_3^2.
\end{equation}
Formula \eqref{eq:moore-two-by-two} is a homogeneous quadratic
polynomial in the six real coordinates
$r,s,c_0,c_1,c_2,c_3$, with coefficients $1$ and $-1$.
\end{example}

\begin{definition}[Quaternionic Monge--Amp\`ere operator]\label{def:qma-operator}
The \emph{quaternionic Monge--Amp\`ere operator} on real $C^2$ functions
is $u\mapsto\det H(u)$.
\end{definition}

By \eqref{eq:trace} and
\eqref{eq:moore-two-by-two},
\[
\begin{aligned}
m=1:&\qquad \det H(u)=\Delta u,\\
m=2:&\qquad \det H(u)
=u_{\bar 1 1}u_{\bar 2 2}-|u_{\bar 1 2}|^2.
\end{aligned}
\]

\begin{lemma}[Positivity and monotonicity]\label{lem:moore-order}
Let $A$ and $B$ be hyperhermitian matrices.
The \emph{positivity property} of the determinant is
\begin{equation}\label{eq:moore-positivity}
A\ge0\ \Longrightarrow\ \det A\ge0,
\qquad A>0\ \Longrightarrow\ \det A>0.
\end{equation}
The \emph{Sylvester criterion} characterizes positive definiteness
by the positivity of all leading principal minors:
\begin{equation}\label{eq:sylvester-criterion}
A>0\quad\Longleftrightarrow\quad
\det\bigl((A_{ab})_{a,b=1}^{\ell}\bigr)>0
\quad\text{for }1\le\ell\le m.
\end{equation}
The determinant is monotone on nonnegative hyperhermitian matrices:
\begin{equation}\label{eq:moore-monotonicity}
0\le B\le A\quad\Longrightarrow\quad\det B\le\det A,
\end{equation}
with strict inequality if $0\le B<A$.
\end{lemma}

\begin{lemma}[Block determinant formula]\label{lem:moore-block}
For $m\ge2$, the determinant of a block matrix satisfies
\begin{equation}\label{eq:moore-block}
\det\begin{pmatrix}B&\zeta\\\zeta^*&c\end{pmatrix}
=\det B\,\bigl(c-\zeta^*B^{-1}\zeta\bigr),
\qquad B>0,\quad \zeta\in\HH^{m-1},\quad c\in\RR,
\end{equation}
where $B$ is an $(m-1)\times(m-1)$ hyperhermitian matrix.
The matrix on the left is positive definite precisely when
$c-\zeta^*B^{-1}\zeta>0$.
\end{lemma}

\begin{proof}
Set
\[
A=\begin{pmatrix}B&\zeta\\\zeta^*&c\end{pmatrix},
\qquad
C=\begin{pmatrix}I_{m-1}&-B^{-1}\zeta\\0&1\end{pmatrix}.
\]
Since $(B^{-1})^*=(B^*)^{-1}=B^{-1}$, we have
\[
C^*=\begin{pmatrix}I_{m-1}&0\\-\zeta^*B^{-1}&1\end{pmatrix}.
\]
Multiplying $A$ first on the left by $C^*$ gives
\begin{equation}\label{eq:block-left-product}
\begin{aligned}
C^*A
&=\begin{pmatrix}
B&\zeta\\
-\zeta^*B^{-1}B+\zeta^*&c-\zeta^*B^{-1}\zeta
\end{pmatrix}\\
&=\begin{pmatrix}B&\zeta\\0&c-\zeta^*B^{-1}\zeta\end{pmatrix}.
\end{aligned}
\end{equation}
Multiplying $C^*A$ in \eqref{eq:block-left-product} on the right
by $C$ gives
\begin{equation}\label{eq:block-diagonalization}
\begin{aligned}
C^*AC
&=\begin{pmatrix}
B&-BB^{-1}\zeta+\zeta\\
0&c-\zeta^*B^{-1}\zeta
\end{pmatrix}\\
&=\begin{pmatrix}B&0\\0&c-\zeta^*B^{-1}\zeta\end{pmatrix}.
\end{aligned}
\end{equation}
The scalar $s:=c-\zeta^*B^{-1}\zeta$ is real, since
\[
\overline{s}
=c-(\zeta^*B^{-1}\zeta)^*
=c-\zeta^*(B^{-1})^*\zeta
=s.
\]
Moreover, $C$ is invertible, with
\[
C^{-1}=\begin{pmatrix}I_{m-1}&B^{-1}\zeta\\0&1\end{pmatrix}.
\]
The real matrix $C_{\RR}$ is upper block triangular with identity
diagonal blocks, so $\det_{\RR}C_{\RR}=1$. Consequently, \eqref{eq:moore-real-determinant} gives
\[
\bigl(\det(C^*C)\bigr)^4
=\det_{\RR}(C_{\RR}^TC_{\RR})
=\bigl(\det_{\RR}C_{\RR}\bigr)^2=1.
\]
Since $C$ is invertible, $C^*C>0$. The positivity property
\eqref{eq:moore-positivity} then gives $\det(C^*C)=1$.
Applying \eqref{eq:moore-congruence} and the eigenvalue formula
\eqref{eq:moore-spectral} to \eqref{eq:block-diagonalization}, we obtain
\[
\det A=\det(C^*AC)
=\det B\,\bigl(c-\zeta^*B^{-1}\zeta\bigr).
\]
For positivity, use \eqref{eq:positivity-congruence}.
By the block diagonalization \eqref{eq:block-diagonalization},
for $\xi=(\eta,t)\in\HH^{m-1}\times\HH$,

\begin{equation}\label{eq:block-quadratic-form}
\xi^*(C^*AC)\xi=\eta^*B\eta+s|t|^2.
\end{equation}
Since $B>0$ and $s\in\RR$, the quadratic form
\eqref{eq:block-quadratic-form} is positive for every
nonzero $\xi$ if and only if $s>0$. Therefore,
\[
A>0\quad\Longleftrightarrow\quad C^*AC>0
\quad\Longleftrightarrow\quad c-\zeta^*B^{-1}\zeta>0.\qedhere
\]
\end{proof}

\subsection{Mixed discriminants and matrix inequalities}
\begin{definition}[Mixed discriminant]\label{def:mixed-determinant}
Let $A_1,\ldots,A_m$ be hyperhermitian $m\times m$ matrices.
The polynomial
\begin{equation}\label{eq:mixed-determinant-definition}
\det(t_1A_1+\cdots+t_mA_m),
\qquad t_1,\ldots,t_m\in\RR,
\end{equation}
is homogeneous of degree $m$.
The \emph{mixed discriminant} $\det(A_1,\ldots,A_m)$ is defined as
$1/m!$ times the coefficient of $t_1\cdots t_m$ in
\eqref{eq:mixed-determinant-definition}.
\end{definition}

The factor $1/m!$ ensures that repeated arguments recover the Moore
determinant. Indeed,
\[
\det\bigl((t_1+\cdots+t_m)A\bigr)
=(t_1+\cdots+t_m)^m\det A
\]
has coefficient $m!\det A$ in front of $t_1\cdots t_m$, so
$\det(A,\ldots,A)=\det A$.
We abbreviate $r$ repetitions of $A$ by $A[r]$, omitting the argument
when $r=0$. The total number of matrix arguments is always $m$.
For example, $\det(B,A[m-1])=\det(B,A,\ldots,A)$.

\begin{example}[Mixed discriminants in dimension two]\label{ex:mixed-determinant}
When $m=2$, write
\[
A=\begin{pmatrix}a&c\\\bar c&b\end{pmatrix},\qquad
B=\begin{pmatrix}d&e\\\bar e&f\end{pmatrix},
\qquad a,b,d,f\in\RR,\quad c,e\in\HH.
\]
Expansion of \eqref{eq:moore-two-by-two} gives, for real $t,s$,
\[
\begin{aligned}
\det(tA+sB)
&=(ta+sd)(tb+sf)-|tc+se|^2\\
&=t^2\det A+ts\bigl(af+bd-2\operatorname{Re}(c\bar e)\bigr)
  +s^2\det B.
\end{aligned}
\]
By Definition~\ref{def:mixed-determinant},
\[
\det(A,B)=\tfrac12\bigl(af+bd-2\operatorname{Re}(c\bar e)\bigr).
\]
\end{example}

\begin{lemma}[Properties of mixed discriminants]\label{lem:mixed-positivity}
The mixed discriminant is symmetric and real-multilinear in its
hyperhermitian arguments, and satisfies
\[
\det(A_1,\ldots,A_m)\ge0\quad\text{if all }A_a\ge0,
\]
with strict positivity if all $A_a>0$. For hyperhermitian $A,B$,
\begin{equation}\label{eq:moore-polarization}
\det(A+sB)=\sum_{r=0}^m\binom mr s^r\det(B[r],A[m-r]),
\qquad s\in\RR.
\end{equation}
For every $C\in\HH^{m\times m}$,
\begin{equation}\label{eq:mixed-congruence}
\det(C^*A_1C,\ldots,C^*A_mC)
=\det(C^*C)\det(A_1,\ldots,A_m).
\end{equation}
Mixed discriminants are monotone in each argument
when the other arguments are nonnegative.
\end{lemma}

\begin{proof}
Permuting the matrices in Definition~\ref{def:mixed-determinant}
permutes the real variables, so the mixed discriminant is symmetric.
For hyperhermitian $A,B$ and $\alpha,\beta\in\RR$, the coefficient
definition gives
\begin{equation}\label{eq:mixed-real-linearity}
\begin{aligned}
\det(\alpha A+\beta B,A_2,\ldots,A_m)
&=\alpha\det(A,A_2,\ldots,A_m)\\
&\quad+\beta\det(B,A_2,\ldots,A_m).
\end{aligned}
\end{equation}
Symmetry gives real linearity in every argument.

For $S\subseteq\{1,\ldots,m\}$, let $|S|$ denote the number of elements
of $S$, and set
\[
C_j^S=\begin{cases}B,&j\in S,\\ A,&j\notin S.\end{cases}
\]
Multilinearity and symmetry give
\begin{equation}\label{eq:mixed-binomial-expansion}
\begin{aligned}
\det(A+sB)
&=\det(A+sB,\ldots,A+sB)\\
&=\sum_{S\subseteq\{1,\ldots,m\}}s^{|S|}\det(C_1^S,\ldots,C_m^S)\\
&=\sum_{r=0}^m s^r
  \sum_{\substack{S\subseteq\{1,\ldots,m\}\\ |S|=r}}
  \det(B[r],A[m-r])\\
&=\sum_{r=0}^m\binom mr s^r\det(B[r],A[m-r]).
\end{aligned}
\end{equation}
Here $\binom mr=m!/(r!(m-r)!)$ counts the subsets with $r$ elements.
This proves \eqref{eq:moore-polarization}.
For \eqref{eq:mixed-congruence}, apply
\eqref{eq:moore-congruence} to $\sum_j t_jA_j$:
\[
\det\left(\sum_jt_jC^*A_jC\right)
=\det(C^*C)\det\left(\sum_jt_jA_j\right),
\]
and compare the coefficients of $t_1\cdots t_m$.

For $A_j\ge0$, choose a unitary matrix
$U_j=(v_{j1},\ldots,v_{jm})$ as in
\eqref{eq:spectral-diagonalization}. Then
\begin{equation}\label{eq:mixed-positive-decomposition}
\begin{aligned}
A_j&=U_j\operatorname{diag}(\lambda_{j1},\ldots,\lambda_{jm})U_j^*\\
&=\sum_{k=1}^m\lambda_{jk}v_{jk}v_{jk}^*,
\qquad \lambda_{jk}\ge0.
\end{aligned}
\end{equation}
By multilinearity, it suffices to prove nonnegativity when
$A_j=v_jv_j^*$, with $v_j\in\HH^m$.
Set $V=(v_1,\ldots,v_m)$. By \eqref{eq:moore-congruence},
\[
\begin{aligned}
\det\left(\sum_{j=1}^m t_jv_jv_j^*\right)
&=\det\bigl(V\operatorname{diag}(t_1,\ldots,t_m)V^*\bigr)\\
&=\det(VV^*)\,t_1\cdots t_m.
\end{aligned}
\]
Comparing coefficients and using $VV^*\ge0$ and
\eqref{eq:moore-positivity} gives
\begin{equation}\label{eq:mixed-rank-one}
\det(v_1v_1^*,\ldots,v_mv_m^*)=\frac{\det(VV^*)}{m!}\ge0.
\end{equation}
Consequently, if $B_1-A_1\ge0$ and $A_2,\ldots,A_m\ge0$, then
\begin{equation}\label{eq:mixed-monotonicity}
\begin{aligned}
\det(B_1,A_2,\ldots,A_m)-\det(A_1,A_2,\ldots,A_m)
&=\det(B_1-A_1,A_2,\ldots,A_m)\\
&\ge0.
\end{aligned}
\end{equation}
Inequality \eqref{eq:mixed-monotonicity} proves monotonicity in the
first argument. Symmetry gives monotonicity in every argument.
If every $A_j>0$, choose $\varepsilon>0$ such that
$A_j\ge\varepsilon I_m$ for all $j$. Monotonicity gives
\[
\det(A_1,\ldots,A_m)
\ge\det(\varepsilon I_m,\ldots,\varepsilon I_m)
=\varepsilon^m\det I_m=\varepsilon^m>0.
\]
\end{proof}

The following identities relate mixed discriminants to the real trace. See
\cite[Proposition 3.2 and Claim 4.7]{Alesker2013}.
We derive the formulas together below.

\begin{lemma}[Variations and concavity]\label{lem:moore-variations}
For $A>0$ and hyperhermitian $B$, the first variation is
\begin{equation}\label{eq:moore-first-variation}
\left.\frac{d}{ds}\right|_{s=0}\det(A+sB)
=m\det(B,A[m-1])
=\det A\,\operatorname{Re}\tr(A^{-1}B).
\end{equation}
For $m\ge2$, the second variation is
\begin{equation}\label{eq:moore-second-variation}
\left.\frac{d^2}{ds^2}\right|_{s=0}\det(A+sB)
=m(m-1)\det(B,B,A[m-2]).
\end{equation}
Put $K=A^{-1/2}BA^{-1/2}$. Then
\begin{equation}\label{eq:moore-log-variations}
\begin{aligned}
\left.\frac{d}{ds}\right|_{s=0}\log\det(A+sB)
&=\operatorname{Re}\tr(A^{-1}B),\\
\left.\frac{d^2}{ds^2}\right|_{s=0}\log\det(A+sB)
&=-\operatorname{Re}\tr(A^{-1}BA^{-1}B)
 =-\tr(K^2)\le0.
\end{aligned}
\end{equation}
Thus $\log\det$ is concave on the positive definite cone.
\end{lemma}

\begin{proof}
Let $\mu_1,\ldots,\mu_m\in\RR$ be the eigenvalues of
$K=A^{-1/2}BA^{-1/2}$. Since
\[
A+sB=A^{1/2}(I_m+sK)A^{1/2},
\]
\eqref{eq:spectral-inverse-square-root}, \eqref{eq:moore-congruence},
and \eqref{eq:spectral-diagonalization} give
\begin{equation}\label{eq:determinant-eigenvalue-expansion}
\det(A+sB)=\det A\,\det(I_m+sK)
=\det A\prod_{a=1}^m(1+s\mu_a).
\end{equation}
For $s$ sufficiently close to zero, every factor in
\eqref{eq:determinant-eigenvalue-expansion} is positive, so
\[
\log\det(A+sB)=\log\det A+\sum_a\log(1+s\mu_a).
\]
Differentiation at zero yields $\sum_a\mu_a=\tr K$ and
$-\sum_a\mu_a^2=-\tr(K^2)$. The cyclic identity \eqref{eq:real-trace-cyclic} gives
\[
\tr K=\operatorname{Re}\tr(A^{-1}B),\qquad
\tr(K^2)=\operatorname{Re}\tr(A^{-1}BA^{-1}B).
\]
The first two derivatives of \eqref{eq:moore-polarization} give the
expressions involving mixed discriminants in
\eqref{eq:moore-first-variation}--\eqref{eq:moore-second-variation}.
The nonpositive second derivative on every line segment contained in
the positive definite cone proves concavity.
\end{proof}
In particular, at $A=I_m$,
\[
m\det(B,I_m[m-1])=\tr B.
\]
For a positive quaternionic Hessian, the first variation reads
\begin{equation}\label{eq:qma-first-variation}
\left.\frac{d}{ds}\right|_{s=0}\log\det H(u+sv)
=\frac{m\det(H(v),H(u)[m-1])}{\det H(u)}
=\operatorname{Re}\tr\bigl(H(u)^{-1}H(v)\bigr).
\end{equation}
Formula \eqref{eq:qma-first-variation} gives the linearized operator
\eqref{eq:linearization} used in the estimates.

\begin{lemma}[Aleksandrov inequality]
\label{lem:moore-aleksandrov}
Let $m\ge2$, let $A>0$ be hyperhermitian, and let $B=B^*$.
Then
\begin{equation}\label{eq:moore-aleksandrov}
\det(B,B,A[m-2])\le
\frac{\det(B,A[m-1])^2}{\det A}.
\end{equation}
Here $B$ need not be nonnegative.
\end{lemma}

\begin{proof}
For completeness, we prove \eqref{eq:moore-aleksandrov} using
\eqref{eq:determinant-eigenvalue-expansion}. Let $\mu_1,\ldots,\mu_m$
be the eigenvalues of $K=A^{-1/2}BA^{-1/2}$. Differentiation gives
\[
\begin{aligned}
\det(B,A[m-1])&=\frac{\det A}{m}\sum_a\mu_a,\\
\det(B,B,A[m-2])
&=\frac{\det A}{m(m-1)}
\left[\left(\sum_a\mu_a\right)^2-\sum_a\mu_a^2\right].
\end{aligned}
\]
The real Cauchy--Schwarz inequality
$(\sum_a\mu_a)^2\le m\sum_a\mu_a^2$ then gives
\eqref{eq:moore-aleksandrov}.
\end{proof}

The rank-one identity below is recorded in
\cite[Claim~4.7]{Alesker2013}.

\begin{lemma}[Rank-one identity and Young inequality]\label{lem:moore-rank-one}
For columns $\xi,\eta\in\HH^m$ and $A>0$,
\begin{equation}\label{eq:moore-rank-one}
\frac{m}{\det A}\det(\xi\xi^*,A[m-1])
=\xi^*A^{-1}\xi\ge0.
\end{equation}
For every $\varepsilon>0$,
\begin{equation}\label{eq:moore-young}
\begin{aligned}
\bigl|\det(\xi\eta^*+\eta\xi^*,A[m-1])\bigr|
&\le\varepsilon\det(\xi\xi^*,A[m-1])\\
&\quad+\varepsilon^{-1}\det(\eta\eta^*,A[m-1]).
\end{aligned}
\end{equation}
\end{lemma}

\begin{proof}
The first variation \eqref{eq:moore-first-variation} and the
real-trace identity \eqref{eq:real-trace-cyclic} give
\[
\frac{m\det(\xi\xi^*,A[m-1])}{\det A}
=\operatorname{Re}\tr(A^{-1}\xi\xi^*)
=\xi^*A^{-1}\xi.
\]
The scalar $\xi^*A^{-1}\xi$ is real and nonnegative because $A^{-1}>0$.
For the second assertion, put
$v_\pm=\sqrt\varepsilon\,\xi\pm\varepsilon^{-1/2}\eta$.
Then
\[
0\le v_\pm v_\pm^*
=\varepsilon\xi\xi^*+\varepsilon^{-1}\eta\eta^*
\pm(\xi\eta^*+\eta\xi^*).
\]
Positivity and multilinearity of the mixed discriminant imply
\eqref{eq:moore-young}.
\end{proof}

\subsection{Affine changes of quaternionic coordinates}\label{sec:quaternionic-coordinates}

\begin{proposition}[Affine covariance]\label{prop:affine-covariance}
Let $U$ be a quaternionic unitary matrix, $U^*U=I_m$, and let
$q_0\in\HH^m$. For a real $C^2$ function $u$, under the affine change
$q=q_0+U\widetilde q$, define
$\widetilde u(\widetilde q)=u(q_0+U\widetilde q)$.
Then
\begin{equation}\label{eq:covariance}
\begin{aligned}
\widetilde u_{\bar a b}(\widetilde q)
&=\sum_{c,d=1}^m\overline{U_{ca}}\,
u_{\bar c d}(q_0+U\widetilde q)U_{db},\\
H(\widetilde u)(\widetilde q)
&=U^*H(u)(q_0+U\widetilde q)U.
\end{aligned}
\end{equation}
The coordinate change preserves Hessian positivity and the
Moore determinant at corresponding points.
\end{proposition}

\begin{proof}
The Hessian transformation can be seen directly from the line-restriction
formula \eqref{eq:compression}. At $q=q_0+U\widetilde q$, for any
$\xi\in\HH^m$,
\begin{equation}\label{eq:hessian-quadratic-covariance}
\xi^*H(\widetilde u)(\widetilde q)\xi
=\left.\Delta_t\,\widetilde u(\widetilde q+\xi t)\right|_{t=0}
=(U\xi)^*H(u)(q)(U\xi)
=\xi^*\bigl(U^*H(u)(q)U\bigr)\xi.
\end{equation}
The second equality uses $U(\xi t)=(U\xi)t$ from
\eqref{eq:matrix-action}. By \eqref{eq:quadratic-recovery}, equality
of quadratic forms in \eqref{eq:hessian-quadratic-covariance} implies
the matrix identity \eqref{eq:covariance}.
By \eqref{eq:positivity-congruence} and
\eqref{eq:moore-unitary}, transformation \eqref{eq:covariance}
preserves Hessian positivity and the Moore determinant at
corresponding points. In particular, the equation transforms as
\[
\det H(u)(q)=f(q)
\quad\Longleftrightarrow\quad
\det H(\widetilde u)(\widetilde q)=f(q_0+U\widetilde q).
\]
\end{proof}

\begin{corollary}[Descent to a flat quotient]\label{cor:hessian-descent}
Let $X=\HH^m/\Gamma$ and let $\pi:\HH^m\to X$ be the quotient
map. For $u\in C^2(M)$ on a domain $M\subset X$, the lift
$\widetilde u=u\circ\pi$ satisfies
\begin{equation}\label{eq:periodic-hessian}
\widetilde u(q+\gamma)=\widetilde u(q),\qquad
H(\widetilde u)(q+\gamma)=H(\widetilde u)(q),\qquad \gamma\in\Gamma.
\end{equation}
By \eqref{eq:periodic-hessian}, the Hessian and the Moore determinant
of the Hessian descend to the quotient.
\end{corollary}

\begin{proof}
Translations correspond to $U=I_m$ in
Proposition~\ref{prop:affine-covariance}. Apply
\eqref{eq:covariance} to the periodic lift $\widetilde u$
to obtain \eqref{eq:periodic-hessian}.
\end{proof}

\par}

\section{Lower-order estimates and reduction to the boundary}\label{sec:estimates}

The subsolution is strictly plurisubharmonic on the compact closure, so
\begin{equation}\label{eq:kappa}
H(\underline u)\ge\kappa I_m\qquad\text{for some }\kappa>0.
\end{equation}
Using the subsolution $\underline u$ in place of $|q|^2$
as a global test function on $M$ is the main technical adjustment in
adapting the estimates from domains in $\HH^m$ to the quotient setting.
Constant real vector fields descend to parallel fields on
$X$, and differentiation along the parallel fields preserves
the form of the equation.

We estimate smooth admissible solutions of
\begin{equation}\label{eq:path}
\det H(u)=F(t):=(1-t)\det H(\underline u)+tf,
\qquad u|_{\partial M}=\phi,
\qquad 0\le t\le1.
\end{equation}
Here $F(t)$ is a function on $\overline M$, and spatial derivatives are
taken with $t$ fixed. The same subsolution works for every $t$. Write
\[
F_-:=\min_{t,\overline M}F(t)>0,\qquad
F_+:=\max_{t,\overline M}F(t).
\]
Throughout the a priori estimates, $C$ denotes a constant depending only
on $m$, $\kappa$, $F_-$, $F_+$, $\sup_t\|\log F(t)\|_{C^2}$,
$\|\underline u\|_{C^4}$, and the fixed geometry of $M\subset X$,
including the boundary and collar bounds for $M$.
The constants $C$ are independent
of $u$ and $t$, and may change from line to line.

\subsection{Comparison and first derivatives}

The comparison principle of Ma-Zhu in \cite[Theorem 2.21]{Zhu} 
is stated for bounded
open subsets of $\HH^m$. For completeness, we prove the
comparison principle in the
$C^2$ setting on the quotient domain considered here.

\begin{lemma}[Comparison principle]\label{lem:comparison}
Let $D\subset M$ be an open subdomain, and let
$v_1,v_2\in C^2(D)\cap C^0(\overline D)$ be plurisubharmonic.
If
\[
\det H(v_1)\ge\det H(v_2)\quad\text{in }D,
\qquad v_1\le v_2\quad\text{on }\partial D,
\]
then $v_1\le v_2$ in $D$.
\end{lemma}

\begin{proof}
Choose $C_0$ so that $p=\underline u-C_0\le-1$ on $\overline M$.
If $v_1-v_2$ has a positive maximum, then, for small $\varepsilon>0$,
$v_1+\varepsilon p-v_2$ has a positive interior maximum. At that point,
\[
H(v_1)+\varepsilon H(\underline u)\le H(v_2),
\]
by \eqref{eq:compression}. {Since $H(v_1)\ge0$ and
$H(\underline u)\ge\kappa I_m>0$, the monotonicity in
Lemma~\ref{lem:moore-order} gives}
\[
\det H(v_2)
\ge\det\bigl(H(v_1)+\varepsilon H(\underline u)\bigr)
>\det H(v_1)
\ge\det H(v_2),
\]
a contradiction.
\end{proof}

Since $\overline M$ is a compact smooth Riemannian manifold with nonempty
boundary, the Dirichlet problem $\Delta h=0$ in $M$, $h=\phi$ on
$\partial M$, has a unique solution $h\in C^\infty(\overline M)$.
Comparison and $\Delta u=\tr H(u)>0$ give
\begin{equation}\label{eq:C0}
\underline u\le u\le h\quad\text{on }\overline M.
\end{equation}
Since $\underline u=u=h=\phi$ on $\partial M$,
\eqref{eq:C0} gives the following inequalities for the inward normal
derivatives:
\[
\partial_\nu\underline u\le\partial_\nu u\le\partial_\nu h.
\]
Tangential derivatives are prescribed by $\phi$. Hence
\begin{equation}\label{eq:boundary-C1}
\|u\|_{C^0(\overline M)}+\sup_{\partial M}|\nabla u|\le C.
\end{equation}

\begin{lemma}\label{lem:first-derivatives}
Let $D$ denote differentiation along the vector field on $X$ induced by
a constant real unit vector field on $\HH^m$.
Then every smooth admissible solution of \eqref{eq:path} satisfies
\[
\sup_M|Du|\le\sup_{\partial M}|Du|+C.
\]
\end{lemma}

\begin{proof}
Let $L$ be the linearization of $v\mapsto\log\det H(v)$ at $u$.
{
By \eqref{eq:qma-first-variation},
\begin{equation}\label{eq:linearization}
Lv=\frac{m}{F(t)}\det(H(v),H(u)[m-1])
  =\operatorname{Re}\sum_{a,b=1}^m u^{a\bar b}v_{\bar b a},
\end{equation}
\par}
where $(u^{a\bar b})_{a,b=1}^m=H(u)^{-1}$.
Since $H(u)>0$, the operator $L$ is elliptic and has no zeroth-order term.
The maximum principle therefore applies. Using
$H(\underline u)\ge\kappa I_m$ and applying the arithmetic--geometric mean
inequality to the eigenvalues of $H(u)^{-1}$, we obtain
\begin{equation}\label{eq:basic-L}
\begin{aligned}
Lu&=m,\qquad
L\underline u=\frac{m}{F(t)}\det(H(\underline u),H(u)[m-1])
             \ge\kappa\sum_{a=1}^m u^{a\bar a},\\
\sum_{a=1}^m u^{a\bar a}&\ge m F(t)^{-1/m}\ge m F_+^{-1/m}>0.
\end{aligned}
\end{equation}
{
Differentiating the logarithm of \eqref{eq:path} along $D$ and using
\eqref{eq:linearization}, we obtain
\begin{equation}\label{eq:differentiated}
L(Du)=D(\log F(t)).
\end{equation}
\par}
Consequently, for a fixed sufficiently large constant $A$,
\[
L(\pm Du+A\underline u)
\ge-C+A\kappa\sum_{a=1}^m u^{a\bar a}\ge0.
\]
The maximum principle gives
$\sup_M|Du|\le\sup_{\partial M}|Du|+A\osc_M\underline u$.
\end{proof}

Combining the lemma with \eqref{eq:boundary-C1} yields
\begin{equation}\label{eq:C1}
\|u\|_{C^1(\overline M)}\le C.
\end{equation}

\subsection{Reduction of the second derivative estimate to the boundary}

\begin{lemma}\label{lem:second-derivatives}
For every derivative $D$ as in Lemma~\ref{lem:first-derivatives},
\begin{equation}\label{eq:interior-C2}
\sup_M D^2u\le\sup_{\partial M}D^2u+C.
\end{equation}
\end{lemma}

\begin{proof}
{
Differentiating the logarithm of \eqref{eq:path} twice along $D$ and
using \eqref{eq:moore-log-variations}, we obtain
\begin{equation}\label{eq:second-differentiated}
\begin{aligned}
L(D^2u)
&=D^2(\log F(t))
 +\operatorname{Re}\tr\bigl(H(u)^{-1}H(Du)H(u)^{-1}H(Du)\bigr)\\
&\ge D^2(\log F(t)).
\end{aligned}
\end{equation}
\par}
By \eqref{eq:basic-L}, choosing $A\ge CF_+^{1/m}/(m\kappa)$ gives
\[
L(D^2u+A\underline u)
\ge-C+A\kappa\sum_{a=1}^m u^{a\bar a}
\ge-C+A\kappa mF_+^{-1/m}\ge0.
\]
The maximum principle proves \eqref{eq:interior-C2}.
\end{proof}

An upper bound for $D^2u$ also gives a lower bound because
$\Delta u>0$. The remaining task is to estimate $D^2u$ on the boundary.

\section{Boundary second derivative estimates}\label{sec:normal}

All estimates in this section concern the quotient domain $M\subset X$.
Let $\pi:\HH^m\to X$ be the quotient map. For each $p\in\partial M$
and $\widetilde p\in\pi^{-1}(p)$, there are open neighborhoods
$U\ni p$ and $\widetilde U\ni\widetilde p$ such that
\[
\pi|_{\widetilde U}:\widetilde U\longrightarrow U
\quad\text{is a diffeomorphism preserving the flat metric},
\qquad
(\pi|_{\widetilde U})^{-1}:U\longrightarrow\widetilde U.
\]
We perform the local computations on
$(\pi|_{\widetilde U})^{-1}(M\cap U)\subset\HH^m$, using the same
notation for functions and their pullbacks. By \eqref{eq:covariance}, quaternionic unitary
affine changes of coordinates preserve the equation. Each chart covers
only a boundary neighborhood. Compactness of $\partial M$ allows a finite
cover and uniform constants in the estimates below.
The local boundary derivative estimates of Ma-Zhu \cite[Sections 3--4]{Zhu}
therefore apply directly. For completeness, we give the detailed proofs below.

The auxiliary functions used below to estimate second derivatives at the
boundary are inspired by the work of B.~Guan~\cite{Guan2009}.

\subsection{A boundary barrier}

For $p\in\overline M$, define
\[
d(p)=\operatorname{dist}_X(p,\partial M).
\]
Choose $\delta_0>0$ sufficiently small that $d$ is smooth up to the boundary
on the collar $\{p\in\overline M:0\le d(p)\le\delta_0\}$, with bounded
second derivatives. There $|\nabla d|=1$, $d=0$ on $\partial M$, and
$\nabla d=\nu$ on $\partial M$, where $\nu$ is the inward unit normal.
All differentiations of $d$ below are restricted to the
collar $\{0\le d\le\delta_0\}$.
Put $v=u-\underline u$. Thus $v\ge0$ and $v=0$ on $\partial M$.

\begin{lemma}\label{lem:barrier}
There are positive constants $N,\eta,\delta$, with $\delta<\delta_0$, such that
\begin{equation}\label{eq:barrier}
V=v+\eta d-\frac N2d^2\ge0,\qquad
LV\le-1-\frac\kappa4\sum_{a=1}^m u^{a\bar a}
\quad\text{for }0\le d\le\delta.
\end{equation}
Moreover, $V=0$ and $0\le\partial_\nu V\le C$ on $\partial M$.
\end{lemma}

\begin{proof}
{
By \eqref{eq:product}, \eqref{eq:linearization}, and
\eqref{eq:moore-rank-one},
\[
L(d^2)=2dLd+2\sum_{a,b=1}^m d_a u^{a\bar b}d_{\bar b}.
\]
\par}
The boundedness of $H(d)$ on the collar and $H(u)^{-1}>0$
give
\[
|Ld|=\left|\operatorname{Re}\tr\bigl(H(u)^{-1}H(d)\bigr)\right|
\le C_d\sum_{a=1}^m u^{a\bar a}
\quad\text{for }0<d<\delta_0.
\]
Since $v=u-\underline u$ and $Lu=m$, we have
\begin{equation}\label{eq:barrier-LV}
\begin{aligned}
LV
&=L(u-\underline u)+\eta Ld-\frac N2L(d^2)\\
&=m-L\underline u+(\eta-Nd)Ld
  -N\sum_{a,b=1}^m d_a u^{a\bar b}d_{\bar b}\\
&\le m-\kappa\sum_{a=1}^m u^{a\bar a}
  +|\eta-Nd|\,|Ld|
  -N\sum_{a,b=1}^m d_a u^{a\bar b}d_{\bar b}\\
&\le m-\bigl(\kappa-C_d(\eta+Nd)\bigr)\sum_{a=1}^m u^{a\bar a}
  -N\sum_{a,b=1}^m d_a u^{a\bar b}d_{\bar b}\\
&=m-\left(\frac\kappa2-C_d(\eta+Nd)\right)\sum_{a=1}^m u^{a\bar a}
  -\left(\frac\kappa2\sum_{a=1}^m u^{a\bar a}
       +N\sum_{a,b=1}^m d_a u^{a\bar b}d_{\bar b}\right),
\end{aligned}
\end{equation}
where we used \eqref{eq:basic-L} and $|\eta-Nd|\le\eta+Nd$.
To estimate the negative bracket in
\eqref{eq:barrier-LV}, let $\lambda_1\le\cdots\le\lambda_m$
be the eigenvalues of $H(u)$.
The smallest eigenvalue of $H(u)^{-1}$ is $\lambda_m^{-1}$, so
{
\begin{equation}\label{eq:inverse-hessian-order}
H(u)^{-1}\ge\lambda_m^{-1}I_m.
\end{equation}
\par}
Applying \eqref{eq:inverse-hessian-order} to the column
$(d_{\bar a})_{a=1}^m$
and using $d_a=\overline{d_{\bar a}}$, we obtain
\[
\sum_{a,b=1}^m d_a u^{a\bar b}d_{\bar b}
\ge\lambda_m^{-1}\sum_{a=1}^m|d_{\bar a}|^2
=\lambda_m^{-1}|\nabla d|^2
=\lambda_m^{-1}.
\]
Moreover,
\[
\sum_{a=1}^m u^{a\bar a}=\sum_{a=1}^m\lambda_a^{-1},
\qquad \prod_{a=1}^m\lambda_a=\det H(u)=F(t)\le F_+.
\]
Dropping the positive term $(\kappa/2)\lambda_m^{-1}$ and applying the
arithmetic--geometric mean inequality to the $m$ positive numbers
$(\kappa/2)\lambda_1^{-1},\ldots,(\kappa/2)\lambda_{m-1}^{-1},
N\lambda_m^{-1}$, we obtain
\begin{equation}\label{eq:barrier-amgm}
\begin{aligned}
&\frac\kappa2\sum_{a=1}^m u^{a\bar a}
  +N\sum_{a,b=1}^m d_a u^{a\bar b}d_{\bar b}\\
&\qquad\ge\frac\kappa2\sum_{a=1}^m\lambda_a^{-1}+N\lambda_m^{-1}\\
&\qquad\ge\sum_{a=1}^{m-1}\frac\kappa{2\lambda_a}+\frac N{\lambda_m}\\
&\qquad\ge m\left(\prod_{a=1}^{m-1}\frac\kappa{2\lambda_a}
                        \cdot\frac N{\lambda_m}\right)^{1/m}\\
&\qquad=m\left(\frac{(\kappa/2)^{m-1}N}{F(t)}\right)^{1/m}
 \ge m\left(\frac{(\kappa/2)^{m-1}N}{F_+}\right)^{1/m}.
\end{aligned}
\end{equation}
Choose $N$ so that the last expression is at least $m+1$, then choose
$\eta$ and $\delta$ so that
\[
C_d(\eta+N\delta)\le\kappa/4,
\qquad N\delta\le\eta.
\]
For $0\le d\le\delta$, \eqref{eq:barrier-LV} and \eqref{eq:barrier-amgm} now give
\[
\begin{aligned}
LV
&\le m-\frac{3\kappa}{4}\sum_{a=1}^m u^{a\bar a}
        -N\sum_{a,b=1}^m d_a u^{a\bar b}d_{\bar b}\\
&=m-\frac\kappa4\sum_{a=1}^m u^{a\bar a}
  -\left(\frac\kappa2\sum_{a=1}^m u^{a\bar a}
          +N\sum_{a,b=1}^m d_a u^{a\bar b}d_{\bar b}\right)\\
&\le m-\frac\kappa4\sum_{a=1}^m u^{a\bar a}-(m+1)
 =-1-\frac\kappa4\sum_{a=1}^m u^{a\bar a}.
\end{aligned}
\]
Also, $v\ge0$ and $N\delta\le\eta$ imply
\[
V=v+d\left(\eta-\frac N2d\right)
\ge v+\frac\eta2d\ge0.
\]
The estimates for $LV$ and $V$ prove \eqref{eq:barrier}. Finally,
$\partial_\nu V=\partial_\nu v+\eta$ on the boundary, and the claim follows
from \eqref{eq:boundary-C1}.
\end{proof}

Fix a boundary point and choose a quaternionic unitary affine chart centered
there so that the inward normal is $\partial_{x_m}$.
Set $\rho=-d$. At the center,
\begin{equation}\label{eq:normal-chart}
\nabla\rho(0)=-\partial_{x_m},\qquad
\rho_{\bar a}(0)=-\delta_{am}.
\end{equation}
On a fixed smaller chart, $|\rho_{x_m}|\ge1/2$.

\begin{lemma}\label{lem:local-comparison}
Let $D=M\cap B_R(0)$ lie in the collar above, and let
$\Psi\in C^2(D)\cap C^1(\overline D)$ be real-valued. Suppose that
\[
L\Psi\ge-C\sum_{a=1}^m u^{a\bar a},\qquad |\Psi|\le C\quad\text{in }D,
\]
and
\[
\Psi\le C|q|^2\quad\text{on }\partial M\cap B_R(0),
\qquad \Psi(0)=0.
\]
Then $\Psi_{x_m}(0)\le C$.
\end{lemma}

\begin{proof}
Choose $B\ge\max\{C,C/R^2\}$. The boundary assumptions and continuity give
\[
\begin{aligned}
\Psi-B|q|^2&\le(C-B)|q|^2\le0
&&\text{on }\partial D\cap\partial M,\\
\Psi-B|q|^2&\le C-BR^2\le0
&&\text{on }\partial D\cap\partial B_R(0).
\end{aligned}
\]
Since $V\ge0$, we obtain
$\Psi-B|q|^2-AV\le0$ on $\partial D$ for every $A\ge0$.
Using $L|q|^2=8\sum_{a=1}^m u^{a\bar a}$ and Lemma~\ref{lem:barrier},
we obtain in $D$
\[
\begin{aligned}
L(\Psi-B|q|^2-AV)
&=L\Psi-8B\sum_{a=1}^m u^{a\bar a}-A LV\\
&\ge-C\sum_{a=1}^m u^{a\bar a}-8B\sum_{a=1}^m u^{a\bar a}
  +A\left(1+\frac\kappa4\sum_{a=1}^m u^{a\bar a}\right)\\
&=A+\left(\frac{A\kappa}{4}-C-8B\right)
      \sum_{a=1}^m u^{a\bar a}\ge0,
\end{aligned}
\]
where we choose $A\ge4(C+8B)/\kappa$.
The maximum principle implies $\Psi\le B|q|^2+AV$ in $D$.
Since $\Psi(0)=V(0)=0$ and $\nu(0)=\partial_{x_m}$, for sufficiently small
$s>0$ the point $(0,\ldots,0,s)\in\HH^m$ lies in $D$, and
\[
\frac{\Psi(0,\ldots,0,s)-\Psi(0)}{s}
\le Bs+A\frac{V(0,\ldots,0,s)-V(0)}{s}.
\]
Letting $s\downarrow0$ and using Lemma~\ref{lem:barrier}, we obtain
\[
\Psi_{x_m}(0)\le A V_{x_m}(0)=A\partial_\nu V(0)\le C.
\]
\end{proof}

\subsection{Tangential and mixed second derivatives}

\medskip
\noindent\textbf{Step 1. Pure tangential derivatives.}
We work in a local lift of a boundary neighborhood to $\HH^m$, using the
same notation as on $M$. Smoothness of $\partial M$ ensures that $\rho=-d$
is a smooth defining function on the inner collar:
{
\begin{equation}\label{eq:collar-rho-properties}
\rho<0\quad\text{in }M,\qquad \rho=0\quad\text{on }\partial M,
\qquad |\nabla\rho|=1.
\end{equation}
\par}
The properties in \eqref{eq:collar-rho-properties} hold
on the collar. No plurisubharmonicity of
$\rho$ is assumed. To factor $v=u-\underline u$, use local collar
coordinates $(z,r)$, where $z$ parametrizes the boundary and $r=\rho$.
Since $v(z,0)=0$,
\[
v(z,r)=\int_0^r\partial_r v(z,s)\,ds
=r\int_0^1\partial_r v(z,tr)\,dt.
\]
Thus $v=\rho\widetilde\gamma$, where
$\widetilde\gamma(z,r)=\int_0^1\partial_r v(z,tr)\,dt$ is smooth up to
the boundary. Differentiating $v=\rho\widetilde\gamma$, we obtain
\[
\nabla v=\widetilde\gamma\nabla\rho+\rho\nabla\widetilde\gamma.
\]
On $\partial M$, we have
\[
\widetilde\gamma=-\partial_\nu v,
\qquad \nabla v=\widetilde\gamma\nabla\rho.
\]
Since $v\ge0$ in $M$ and $v=0$ on $\partial M$, we have
$\partial_\nu v\ge0$. The boundary gradient estimate
\eqref{eq:boundary-C1} therefore gives $-C\le\widetilde\gamma\le0$.
For $p\in\partial M$ and real unit vectors $e,\tau\in T_p\partial M$,
the product rule gives
\[
\begin{aligned}
D^2v(e,\tau)
&=\widetilde\gamma D^2\rho(e,\tau)
  +(D_e\rho)(D_\tau\widetilde\gamma)
  +(D_\tau\rho)(D_e\widetilde\gamma)
  +\rho D^2\widetilde\gamma(e,\tau).
\end{aligned}
\]
At $p$, $\rho(p)=D_e\rho(p)=D_\tau\rho(p)=0$, so
\begin{equation}\label{eq:tangential}
D^2v(e,\tau)(p)=\widetilde\gamma(p)D^2\rho(e,\tau)(p).
\end{equation}
Consequently,
\[
\begin{aligned}
|D^2u(e,\tau)(p)|
&\le |D^2\underline u(e,\tau)(p)|
   +|\widetilde\gamma(p)|\,|D^2\rho(e,\tau)(p)|\\
&\le\|D^2\underline u\|_{L^\infty(\partial M)}
   +C\|D^2\rho\|_{L^\infty(\partial M)}\le C.
\end{aligned}
\]

\medskip
\noindent\textbf{Step 2. Mixed tangential-normal derivatives.}

\begin{proposition}\label{prop:mixed}
For every real unit tangent vector $\tau$ on $\partial M$,
$|D^2u(\tau,\nu)|\le C$.
\end{proposition}

\begin{proof}
For any real coordinate $r\in\{x_a,y_a,z_a,w_a:1\le a\le m\}$,
the fixed $C^3$ bound for $\underline u$ gives
$-CI_m\le H(\underline u_r)\le CI_m$. Since $v=u-\underline u$,
\eqref{eq:linearization} and \eqref{eq:differentiated} yield
\begin{equation}\label{eq:Lv-r}
\begin{aligned}
|Lv_r|
&=\left|(\log F(t))_r
       -\operatorname{Re}\tr\bigl(H(u)^{-1}H(\underline u_r)\bigr)\right|\\
&\le |(\log F(t))_r|
       +\left|\operatorname{Re}\tr\bigl(H(u)^{-1}H(\underline u_r)\bigr)\right|\\
&\le C+C\sum_{a=1}^m u^{a\bar a}
 \le C\sum_{a=1}^m u^{a\bar a}.
\end{aligned}
\end{equation}
The last inequality uses \eqref{eq:basic-L}, namely
$1\le(F_+^{1/m}/m)\sum_a u^{a\bar a}$.

In the coordinates chosen in \eqref{eq:normal-chart}, let $\tau$ denote any coordinate direction other than
$x_m$. Define
\[
\beta=\frac{\rho_\tau}{\rho_{x_m}},\qquad
T=\partial_\tau-\beta\partial_{x_m}.
\]
Here $\beta$ is a real-valued smooth function with $\|\beta\|_{C^2}\le C$,
and $\beta_a=\partial\beta/\partial q_a$ as in
\eqref{eq:quaternionic-derivatives}.
Since $T\rho=\rho_\tau-\beta\rho_{x_m}=0$, $T$ is tangent to
$\partial M$ and $Tv=0$ there. Also, $\beta(0)=0$, so
\[
(Tv)_{x_m}(0)=v_{\tau x_m}(0)-\beta_{x_m}(0)v_{x_m}(0).
\]
The term $\beta_{x_m}(0)v_{x_m}(0)$ is bounded by
\eqref{eq:boundary-C1}. Therefore, we only need
to bound $|(Tv)_{x_m}(0)|$ by Lemma~\ref{lem:local-comparison}.
{
By \eqref{eq:product}, \eqref{eq:real-trace-cyclic}, and
\eqref{eq:linearization},
\begin{equation}\label{eq:mixed-L}
\begin{aligned}
L(Tv)&=Lv_\tau-\beta Lv_{x_m}-v_{x_m}L\beta\\
&\quad-2\operatorname{Re}\sum_{a,b=1}^m
 \beta_a u^{a\bar b}v_{x_m\bar b}.
\end{aligned}
\end{equation}
\par}
Using \eqref{eq:Lv-r}, \eqref{eq:C1}, and the fixed $C^2$ bound for
$\beta$, we obtain
\begin{equation}\label{eq:LTv-bound}
|L(Tv)|\le C\sum_{a=1}^m u^{a\bar a}
 +2\left|\operatorname{Re}\sum_{a,b=1}^m
       \beta_a u^{a\bar b}v_{x_m\bar b}\right|.
\end{equation}
Applying $\partial/\partial\bar q_b$ to
$v_m=v_{x_m}-v_{y_m}i-v_{z_m}j-v_{w_m}k$ gives
\begin{equation}\label{eq:column}
v_{x_m\bar b}=v_{\bar b m}
                 +v_{y_m\bar b}i+v_{z_m\bar b}j+v_{w_m\bar b}k,
\qquad 1\le b\le m.
\end{equation}
Insert
$v_{\bar b m}=u_{\bar b m}-\underline u_{\bar b m}$ into
\eqref{eq:column}. The inverse-matrix identity gives
{
\begin{equation}\label{eq:normal-cancellation}
\operatorname{Re}\sum_{a,b=1}^m\beta_a u^{a\bar b}u_{\bar b m}
=\operatorname{Re}\beta_m=\beta_{x_m}.
\end{equation}
\par}
The term $\beta_{x_m}$ in \eqref{eq:normal-cancellation}
and the subsolution term are bounded by $C\sum_a u^{a\bar a}$,
using \eqref{eq:basic-L}.
{
For any real coordinate $r$ and $\omega\in\{1,i,j,k\}$, apply
\eqref{eq:moore-rank-one}--\eqref{eq:moore-young} with $A=H(u)$
to the columns $(\beta_{\bar a})_{a=1}^m$ and
$(v_{r\bar a}\omega)_{a=1}^m$.
The rank-one matrix associated with the second column has entries
$v_{r\bar a}\omega\bar\omega v_{rb}=v_{r\bar a}v_{rb}$.
The fixed $C^2$ bound for $\beta$ therefore gives, for every
$\varepsilon>0$,
\begin{equation}\label{eq:cross-young}
\begin{aligned}
\left|\operatorname{Re}\sum_{a,b=1}^m
 \beta_a u^{a\bar b}(v_{r\bar b}\omega)\right|
&\le C_\varepsilon\sum_{a=1}^m u^{a\bar a}\\
&\quad+\frac{\varepsilon m}{F(t)}
 \det((v_{r\bar a}v_{rb})_{a,b},H(u)[m-1]).
\end{aligned}
\end{equation}
Combining \eqref{eq:cross-young} with \eqref{eq:column} and
\eqref{eq:normal-cancellation}, and rescaling $\varepsilon$, yields
\par}
\begin{equation}\label{eq:cross}
\begin{aligned}
\left|\operatorname{Re}\sum_{a,b=1}^m\beta_a u^{a\bar b}v_{x_m\bar b}\right|
&\le C_\varepsilon\sum_{a=1}^m u^{a\bar a}
  +\frac{\varepsilon m}{F(t)}\sum_{r\in\{y_m,z_m,w_m\}}
        \det((v_{r\bar a}v_{rb})_{a,b},H(u)[m-1]),\\
\left|\operatorname{Re}\sum_{a,b=1}^m\beta_a u^{a\bar b}v_{r\bar b}\right|
&\le C_\varepsilon\sum_{a=1}^m u^{a\bar a}
  +\frac{\varepsilon m}{F(t)}
        \det((v_{r\bar a}v_{rb})_{a,b},H(u)[m-1]),
\quad r\in\{y_m,z_m,w_m\}.
\end{aligned}
\end{equation}

The estimates \eqref{eq:cross} hold for any real-valued smooth
$\beta$ with a fixed $C^2$ bound, and also with the coordinate index
$m$ replaced by any $1\le \ell\le m$.
{
By \eqref{eq:product}, \eqref{eq:linearization}, \eqref{eq:Lv-r}, and
\eqref{eq:C1},
\begin{equation}\label{eq:derivative-squares}
\begin{aligned}
L(v_r^2)
&=2v_rLv_r+\frac{2m}{F(t)}
       \det((v_{r\bar a}v_{rb})_{a,b},H(u)[m-1])\\
&\ge\frac{2m}{F(t)}\det((v_{r\bar a}v_{rb})_{a,b},H(u)[m-1])
       -C\sum_{a=1}^m u^{a\bar a}.
\end{aligned}
\end{equation}
\par}
We therefore define
\[
\Psi=\pm Tv+v_{y_m}^2+v_{z_m}^2+v_{w_m}^2.
\]
Combining \eqref{eq:LTv-bound}, \eqref{eq:cross}, and
\eqref{eq:derivative-squares}, we obtain
\[
\begin{aligned}
L\Psi
&\ge-|L(Tv)|+\sum_{r\in\{y_m,z_m,w_m\}}L(v_r^2)\\
&\ge-C_\varepsilon\sum_{a=1}^m u^{a\bar a}
 +\frac{(2-2\varepsilon)m}{F(t)}
   \sum_{r\in\{y_m,z_m,w_m\}}
       \det((v_{r\bar a}v_{rb})_{a,b},H(u)[m-1])\\
&\ge-C\sum_{a=1}^m u^{a\bar a},
\end{aligned}
\]
{where we take $\varepsilon=1/2$ and use
Lemma~\ref{lem:mixed-positivity} for the nonnegativity of the mixed discriminants.}
The function $\Psi$ is bounded by \eqref{eq:C1}. On
$\partial M\cap B_R$, we have $Tv=0$ and
$v_r=\widetilde\gamma\rho_r$ for $r\in\{y_m,z_m,w_m\}$.
Since $|\widetilde\gamma|\le C$, $\rho_r(0)=0$, and $\rho$ has a fixed
$C^2$ bound,
\[
|\rho_r(q)|=|\rho_r(q)-\rho_r(0)|\le C|q|,
\qquad |v_r(q)|\le C|q|.
\]
Consequently,
\[
0\le\Psi(q)=\sum_{r\in\{y_m,z_m,w_m\}}v_r(q)^2
\le C|q|^2\quad(q\in\partial M\cap B_R),
\qquad\Psi(0)=0.
\]
At the origin, $v_r(0)=0$ for $r=y_m,z_m,w_m$, and hence
\[
\partial_{x_m}(v_r^2)(0)=2v_r(0)v_{rx_m}(0)=0,
\qquad \Psi_{x_m}(0)=\pm(Tv)_{x_m}(0).
\]
Thus Lemma~\ref{lem:local-comparison}, applied for both signs, gives
\[
\pm(Tv)_{x_m}(0)\le C,
\qquad |v_{\tau x_m}(0)|
\le |(Tv)_{x_m}(0)|+|\beta_{x_m}(0)v_{x_m}(0)|\le C.
\]
Adding the fixed derivative $\underline u_{\tau x_m}(0)$ proves the claim.
\end{proof}

\subsection{The double normal derivative}

The normal estimate requires a positive lower bound for the
quaternionic tangential Hessian. We first prove
\eqref{eq:tangent-positive}, following the
minimizing-direction argument of Ma-Zhu in \cite[Section 4]{Zhu}.

At $p\in\partial M$, the quaternionic tangent space is the orthogonal
complement of $(\rho_{\bar a}(p))_{a=1}^m$ with respect to the quaternionic
Hermitian inner product, since
\[
\bigl((\rho_{\bar a}(p))_{a=1}^m\bigr)^*\xi
=\sum_{a=1}^m\overline{\rho_{\bar a}(p)}\xi_a
=\sum_{a=1}^m\rho_a(p)\xi_a.
\]
Thus
{
\begin{equation}\label{eq:quaternionic-tangent-space}
T_p^{\HH}\partial M
=\bigl((\rho_{\bar a}(p))_{a=1}^m\bigr)^\perp
=\left\{\xi\in\HH^m:\sum_{a=1}^m\rho_a(p)\xi_a=0\right\}.
\end{equation}
\par}
The tangent space \eqref{eq:quaternionic-tangent-space}
has quaternionic dimension $m-1$.

\begin{proposition}\label{prop:tangent-positive}
There is a constant $c_0>0$ such that every solution of
\eqref{eq:path} satisfies
\begin{equation}\label{eq:tangent-positive}
\xi^*H(u)(p)\xi\ge c_0|\xi|^2,
\qquad p\in\partial M,\quad \xi\in T_p^{\HH}\partial M.
\end{equation}
\end{proposition}

\begin{proof}
Let
\[
\mu_0=\min_{\substack{p\in\partial M,\ |\xi|=1\\
                      \xi\in T_p^{\HH}\partial M}}\xi^*H(u)(p)\xi.
\]
We must bound $\mu_0$ below independently of the solution. If
$\mu_0\ge\kappa/2$, there is nothing to prove, so assume
$0<\mu_0<\kappa/2$. Choose a minimizing point as the origin and use
\eqref{eq:normal-chart}, with minimizing direction $(1,0,\ldots,0)^T$.
Thus $\mu_0=u_{\bar1 1}(0)$.
To prove \eqref{eq:tangent-positive}, we only need an
upper bound for $u_{x_mx_m}(0)$ at the minimizing point.
Indeed, if $u_{x_mx_m}(0)\le C$, then \eqref{eq:tangential} gives
\[
\begin{aligned}
0<\tr H(u)(0)=\Delta u(0)
&=u_{x_mx_m}(0)+\sum_{a=1}^{m-1}u_{x_ax_a}(0)\\
&\quad+\sum_{a=1}^m
 \bigl(u_{y_ay_a}(0)+u_{z_az_a}(0)+u_{w_aw_a}(0)\bigr)
\le C.
\end{aligned}
\]
Writing $0<\lambda_1\le\cdots\le\lambda_m$ for the eigenvalues of
$H(u)(0)$, we obtain
\[
\lambda_m\le C,\qquad
\lambda_1=\frac{\det H(u)(0)}{\lambda_2\cdots\lambda_m}
\ge\frac{F_-}{C^{m-1}},\qquad
\mu_0=u_{\bar1 1}(0)\ge\lambda_1.
\]
The lower bound for $\mu_0$ yields
\eqref{eq:tangent-positive}. We now prove $u_{x_mx_m}(0)\le C$
at the minimizing point.

In the local lift, choose $R>0$ so that
$\overline D=\overline{M\cap B_R(0)}$ lies in the coordinate chart and
the collar where $\rho=-d$ is smooth. To use the definition of $\mu_0$
on $\partial M\cap B_R(0)$, we seek a smooth field
$\xi:\overline D\to\HH^m$ satisfying
{
\begin{equation}\label{eq:tangent-field-conditions}
|\xi|=1,\qquad \sum_{a=1}^m\rho_a\xi_a=0,\qquad
\xi_a(0)=\delta_{a1},\qquad \|\xi\|_{C^2(\overline D)}\le C.
\end{equation}
\par}
The conditions \eqref{eq:tangent-field-conditions} and
the definition of $\mu_0$ give $\xi^*H(u)\xi\ge\mu_0$ on
$\partial M\cap B_R(0)$, with equality at the origin.
For each $q\in\overline D$, project $(1,0,\ldots,0)^T$ onto the
quaternionic tangent space of the level set $\{\rho=\rho(q)\}$ at $q$.
Since $\sum_a|\rho_{\bar a}|^2=1$, the orthogonal projection is
\[
\left(I_m-(\rho_{\bar a}\rho_b)_{a,b=1}^m\right)
(1,0,\ldots,0)^T
=(\delta_{a1}-\rho_{\bar a}\rho_1)_{a=1}^m,
\]
and the squared norm of the projected vector is
\[
\sum_{a=1}^m|\delta_{a1}-\rho_{\bar a}\rho_1|^2
=1-|\rho_{\bar1}|^2.
\]
Normalizing defines $\xi$ on all of $\overline D$:
\begin{equation}\label{eq:tangent-field}
\xi_a(q)=\frac{\delta_{a1}-\rho_{\bar a}(q)\rho_1(q)}
                   {\sqrt{1-|\rho_{\bar1}(q)|^2}},
\qquad q\in\overline D,\quad 1\le a\le m.
\end{equation}
Since $\rho_{\bar1}(0)=0$, reducing $R$ makes the denominator at least
$1/2$ throughout $\overline D$. The stated properties follow. In
particular, $\xi(q)\in T_q^{\HH}\partial M$ for
$q\in\partial M\cap\overline D$.

To express the boundary identity for $u-\underline u$, define
\[
\alpha=\xi^*H(\underline u)\xi,\qquad
\beta=\xi^*H(\rho)\xi,\qquad
\gamma=\langle\nabla v,\nabla\rho\rangle.
\]
Here the gradients and inner product are real Euclidean, so $\gamma$ is
real-valued. Equivalently,
\[
\gamma=\operatorname{Re}\sum_{a=1}^m v_a\rho_{\bar a}
=\frac12\sum_{a=1}^m
 \bigl(v_a\rho_{\bar a}+\rho_a v_{\bar a}\bigr).
\]
The factorization $v=\rho\widetilde\gamma$ from Step 1 gives
$\gamma=\widetilde\gamma$ on $\partial M$.
{
For $\omega\in\{1,i,j,k\}$, \eqref{eq:tangent-field-conditions} gives
\[
D\rho[\xi\omega]
=\operatorname{Re}\left[\left(\sum_{a=1}^m\rho_a\xi_a\right)
\omega\right]=0.
\]
Thus $\xi,\xi i,\xi j,\xi k$ are real tangent vectors on the boundary.
Applying \eqref{eq:compression} and \eqref{eq:tangential}, we obtain
on $\partial M\cap B_R$
\begin{equation}\label{eq:boundary-compression}
\xi^*H(u)\xi
=\xi^*H(\underline u)\xi+\gamma\xi^*H(\rho)\xi
=\alpha+\beta\gamma.
\end{equation}
\par}
At the center $\alpha(0)\ge\kappa$ and
$\alpha(0)+\beta(0)\gamma(0)=\mu_0$.
Since $-C_1\le\gamma\le0$ on $\partial M$, with $C_1\ge1$ fixed, it follows that
\begin{equation}\label{eq:b-positive}
\beta(0)\ge\frac{\kappa}{2C_1}>0.
\end{equation}

The minimizing property gives $\alpha+\beta\gamma\ge\mu_0$ on
$\partial M\cap B_R$, with equality at the origin. We therefore define
{
\begin{equation}\label{eq:normal-auxiliary-function}
\Phi=\mu_0-\alpha-\beta\gamma,
\qquad \Phi\le0\ \text{on }\partial M\cap B_R,
\qquad \Phi(0)=0.
\end{equation}
\par}
The function $\Phi$ in \eqref{eq:normal-auxiliary-function}
relates the desired normal derivative to
Lemma~\ref{lem:local-comparison}. Indeed, $|\nabla\rho|=1$ gives
$D^2\rho\,\nabla\rho=0$, and hence $\rho_{rx_m}(0)=0$ for every real
coordinate $r$. Therefore
{
\begin{equation}\label{eq:normal-auxiliary-derivative}
\begin{aligned}
\gamma_{x_m}(0)&=-v_{x_mx_m}(0),\\
\Phi_{x_m}(0)&=-\alpha_{x_m}(0)-\beta_{x_m}(0)\gamma(0)
                    +\beta(0)v_{x_mx_m}(0).
\end{aligned}
\end{equation}
\par}
By \eqref{eq:b-positive} and
\eqref{eq:normal-auxiliary-derivative}, an upper bound for
$v_{x_mx_m}(0)$ follows from
\[
\Phi_{x_m}(0)\le C.
\]
As in the mixed derivative estimate, we add squares of first derivatives
to apply Lemma~\ref{lem:local-comparison}.

We first estimate $L\Phi$. Writing $r_1,\ldots,r_{4m}$ for the real
coordinates, we have
\[
\Phi=\mu_0-\alpha-\sum_{j=1}^{4m}\beta\rho_{r_j}v_{r_j}.
\]
By the product rule,
\[
\begin{aligned}
L\Phi=-L\alpha
&-\sum_{j=1}^{4m}\bigl(\beta\rho_{r_j}Lv_{r_j}
                         +v_{r_j}L(\beta\rho_{r_j})\bigr)\\
&-2\operatorname{Re}\sum_{j=1}^{4m}\sum_{a,b=1}^m
             (\beta\rho_{r_j})_a u^{a\bar b}v_{r_j\bar b}.
\end{aligned}
\]
The coefficients have fixed $C^2$ bounds. Applying \eqref{eq:cross} to
each term in the last sum gives, for every $\varepsilon>0$,
\[
L\Phi\ge-C_\varepsilon\sum_{a=1}^m u^{a\bar a}
-\frac{\varepsilon m}{F(t)}\sum_{k=1}^m\sum_{r\in\{y_k,z_k,w_k\}}
       \det((v_{r\bar a}v_{rb})_{a,b},H(u)[m-1]).
\]
The mixed discriminants are supplied by \eqref{eq:derivative-squares}.
Define
\[
\Psi=\Phi+\sum_{k=1}^m(v_{y_k}^2+v_{z_k}^2+v_{w_k}^2).
\]
Then, taking $\varepsilon=1$ in the preceding estimate,
\[
L\Psi\ge-C\sum_{a=1}^m u^{a\bar a}
 +\frac{m}{F(t)}\sum_{k=1}^m\sum_{r\in\{y_k,z_k,w_k\}}
 \det((v_{r\bar a}v_{rb})_{a,b},H(u)[m-1])
\ge-C\sum_{a=1}^m u^{a\bar a}.
\]
On $\partial M\cap B_R$, $\Phi\le0$ and
$|v_r|=|\gamma\rho_r|\le C|q|$ for $r=y_a,z_a,w_a$,
since $\rho_r(0)=0$ and $\rho$ has a fixed $C^2$ bound.
Consequently,
\[
\Psi\le C|q|^2\quad\text{on }\partial M\cap B_R,
\qquad \Psi(0)=0.
\]
The function $\Psi$ is uniformly bounded by \eqref{eq:C1}. Here
$0<\mu_0<\kappa/2$. Lemma~\ref{lem:local-comparison} implies
$\Phi_{x_m}(0)\le C$, since the added squares have zero normal derivative
at the center.

The normal derivative identity above and \eqref{eq:b-positive} now give
$u_{x_mx_m}(0)\le C$, as required.
\end{proof}

\begin{proposition}\label{prop:C2}
Every smooth admissible solution of \eqref{eq:path} satisfies
\begin{equation}\label{eq:C2}
\|u\|_{C^2(\overline M)}\le C,\qquad
\lambda_-I_m\le H(u)\le\lambda_+I_m
\end{equation}
for fixed positive constants $\lambda_-,\lambda_+$.
\end{proposition}

\begin{proof}
We bound the double normal derivative at every boundary
point using the determinant argument of Ma-Zhu \cite[Lemma 4.2]{Zhu}.
Choose the coordinates of \eqref{eq:normal-chart} and diagonalize the
quaternionic tangential block by a unitary change in the first $m-1$
coordinates. At the center,
\[
u_{\bar a b}=\lambda_a\delta_{ab}\quad(1\le a,b\le m-1),
\qquad c_0\le\lambda_a\le C,
\qquad |u_{\bar a m}|\le C\quad(1\le a\le m-1),
\]
by Proposition~\ref{prop:tangent-positive}, \eqref{eq:tangential}, and
Proposition~\ref{prop:mixed}.
{The determinant formula for block matrices in
Lemma~\ref{lem:moore-block} gives}
\begin{equation}\label{eq:schur}
\begin{aligned}
F(t)&=\left(\prod_{a=1}^{m-1}\lambda_a\right)
 \left(u_{\bar m m}-\sum_{a=1}^{m-1}
          \frac{|u_{\bar a m}|^2}{\lambda_a}\right),\\
0<u_{\bar m m}
&=\frac{F(t)}{\prod_{a=1}^{m-1}\lambda_a}
  +\sum_{a=1}^{m-1}\frac{|u_{\bar a m}|^2}{\lambda_a}\le C.
\end{aligned}
\end{equation}
Since $u_{y_my_m}$, $u_{z_mz_m}$, and $u_{w_mw_m}$ are bounded,
\[
u_{x_mx_m}=u_{\bar m m}-u_{y_my_m}-u_{z_mz_m}-u_{w_mw_m},
\qquad |u_{x_mx_m}|\le C.
\]

All boundary second derivatives are now bounded. By
\eqref{eq:interior-C2}, $D^2u\le CI_{4m}$ throughout $M$.
Since $\Delta u>0$, every real eigenvalue is also at least $-(4m-1)C$.
The upper and lower bounds for $D^2u$ prove the $C^2$
estimate in \eqref{eq:C2}. By \eqref{eq:compression},
$0<H(u)\le4CI_m$. The determinant lower bound then gives
\[
\lambda_{\min}(H(u))\ge F_-(4C)^{-(m-1)}.
\]
\end{proof}

\section{Proofs of Theorem~\ref{thm:main} and Corollary~\ref{cor:three}}\label{sec:existence}

{The operator $u\mapsto\log\det H(u)$ is concave by
Lemma~\ref{lem:moore-variations}. Estimate \eqref{eq:C2} bounds the
eigenvalues of $H(u)$ above and away from zero, so the
linearized operator $L$ in \eqref{eq:linearization} is uniformly elliptic, with constants independent of $t$.}
Indeed, for every quadratic polynomial
$h$ with $D^2h\ge0$, \eqref{eq:compression} and \eqref{eq:trace} give
\[
\lambda_+^{-1}\tr D^2h
\le Lh=\operatorname{Re}\tr\bigl(H(u)^{-1}H(h)\bigr)
\le\lambda_-^{-1}\tr D^2h.
\]
To apply the standard regularity estimates, we extend $\log\det$
outside the eigenvalue range in \eqref{eq:C2} without changing the
equation satisfied by $u$. Choose a smooth concave function
$\chi:\RR\to\RR$ such that
\[
\chi(s)=\log s\quad\text{near }[\lambda_-,\lambda_+],
\qquad 0<c\le\chi'(s)\le C\quad(s\in\RR).
\]
For any hyperhermitian matrix $A$ with eigenvalues
$\lambda_1(A),\ldots,\lambda_m(A)$, set
{
\begin{equation}\label{eq:elliptic-extension}
G(A)=\sum_{a=1}^m\chi(\lambda_a(A)).
\end{equation}
\par}
The function $G$ in \eqref{eq:elliptic-extension} is
smooth and concave. For every hyperhermitian $B\ge0$,
\[
c\tr B\le\left.\frac{d}{d\varepsilon}\right|_{\varepsilon=0}
G(A+\varepsilon B)\le C\tr B.
\]
Since $H$ depends linearly on the real Hessian, and
$D^2h\ge0$ implies $H(h)\ge0$ with $\tr H(h)=\Delta h$,
the operator $v\mapsto G(H(v))$ is smooth, concave, and uniformly
elliptic in the real Hessian. By \eqref{eq:C2}, every solution of
\eqref{eq:path} satisfies
\[
G(H(u))=\log\det H(u)=\log F(t).
\]
Applying the Evans--Krylov interior estimate
\cite{Evans,Krylov1} and the boundary Dirichlet estimate
\cite{Krylov2} to the uniformly elliptic concave equation
$G(H(u))=\log F(t)$ in local charts yields, for some
$\alpha\in(0,1)$,
\begin{equation}\label{eq:holder}
\|u\|_{C^{2,\alpha}(\overline M)}\le C.
\end{equation}
Linear Schauder estimates \cite[Chapters 6 and 17]{GT} then give uniform
higher-order bounds by differentiating the equation and using the smooth
boundary data.

\begin{proof}[Proof of Theorem~\ref{thm:main}]
The path \eqref{eq:path} starts at $u=\underline u$.
On the open set of $w\in C^{2,\alpha}(\overline M)$ with zero boundary
trace and $H(\underline u+w)>0$, consider
{
\begin{equation}\label{eq:continuity-map}
\mathcal N(t,w)=\log\det H(\underline u+w)-\log F(t).
\end{equation}
\par}
The derivative $D_w\mathcal N$ of
\eqref{eq:continuity-map} is the operator $L$ in
\eqref{eq:linearization}. At a smooth admissible solution, $L$ is a
uniformly elliptic operator with smooth coefficients and no zeroth-order
term. The maximum principle and linear Dirichlet theory imply that
\[
L:C^{2,\alpha}_0(\overline M)\longrightarrow C^\alpha(\overline M)
\]
is an isomorphism. The subscript denotes zero boundary trace.
The implicit function theorem proves openness of the set of solvable
parameters. Estimates \eqref{eq:C2} and \eqref{eq:holder}, followed by
compactness in $C^{2,\beta}$ for $0<\beta<\alpha$, prove closedness.
The positive lower bound for $H(u)$ persists in the limit. Regularity
therefore gives a smooth solution at $t=1$. Uniqueness follows from
Lemma~\ref{lem:comparison}.
\end{proof}

\begin{remark}
Under the hypotheses of Theorem~\ref{thm:appendix},
a smooth strictly plurisubharmonic subsolution
$\underline u$ exists automatically.
Replacing $e^{\lambda|q|^2}$ in the auxiliary functions of Ma-Zhu 
\cite[Lemmas~3.2 and~3.3]{Zhu} with $A\underline u$,
for sufficiently large $A$, yields the smooth Dirichlet
solvability in Theorem~\ref{thm:appendix}.
For completeness, we have proved the general existence result
under a subsolution hypothesis.
\end{remark}

\begin{proof}[Proof of Corollary~\ref{cor:three}]
Add a torus in the $w$ coordinates and set
\[
\widetilde X=(\RR^m/\Lambda)\times\RR^{2m}\times(\RR^m/\ZZ^m),
\qquad \widetilde\Omega=\Omega\times\TT^m.
\]
The space $\widetilde X$ is a flat quaternionic quotient
with translation group
$\{\lambda+kn:\lambda\in\Lambda,\ n\in\ZZ^m\}$.
Extend the data and the subsolution to
$\widetilde\Omega$ as functions independent of $w$.
Then $H(\underline u)=H_3(\underline u)$, so Theorem~\ref{thm:main}
gives a unique smooth solution on $\widetilde\Omega$.
Since the lifted domain and data are invariant under all
translations in the added torus, uniqueness implies that the solution is
independent of $w$ and descends to $\Omega$.
Comparison and Remark~\ref{prop:families} descend in the same way.
\end{proof}

\begin{remark}\label{prop:families}
Let $M_\theta\Subset X$ and the data satisfy the hypotheses of
Theorem~\ref{thm:main}. Suppose that
$T_\theta:\overline M_0\to\overline M_\theta$ is a smooth family of
diffeomorphisms, with $T_0=\mathrm{id}$, and that the pullbacks
$f_\theta\circ T_\theta$ and $\underline u_\theta\circ T_\theta$ are smooth
in $(\theta,q)$ up to the boundary, where $\theta$ varies in a
finite-dimensional parameter space. Set
$w_\theta=(u_\theta-\underline u_\theta)\circ T_\theta$.
On the fixed domain $M_0$, the equation becomes
{
\begin{equation}\label{eq:family-nonlinear-map}
\mathcal N(\theta,w)=
\left[\log\det H\bigl(\underline u_\theta+w\circ T_\theta^{-1}\bigr)
      -\log f_\theta\right]\circ T_\theta=0,
\qquad w|_{\partial M_0}=0.
\end{equation}
\par}
The map $\mathcal N$ in \eqref{eq:family-nonlinear-map}
is smooth where the Hessian is positive, and
{
\begin{equation}\label{eq:family-linearization}
D_w\mathcal N(0,w_0)[h]
=\operatorname{Re}\tr\bigl(H(u_0)^{-1}H(h)\bigr).
\end{equation}
\par}
As in the proof of Theorem~\ref{thm:main}, the derivative
\eqref{eq:family-linearization} is an
isomorphism from $C^{k+2,\alpha}_0(\overline M_0)$ to
$C^{k,\alpha}(\overline M_0)$ for every $k\ge0$ and $0<\alpha<1$.
The same implicit function argument and uniqueness therefore show that
\[
(\theta,q)\longmapsto u_\theta(T_\theta(q))
\]
is smooth up to the boundary, with locally uniform solution norms.
\end{remark}

\section{ \cite[Theorem A.1]{OpenAI} as a special case of  Corollary~\ref{cor:three}}\label{sec:descent}

Take the data in Theorem~\ref{thm:appendix}.
For sufficiently large $b$, $B_0+br$ is admissible and
$\det H_3(B_0+br)\ge f$.
Set $u=B_0+v$. For $m=1$, the equation reduces to
\[
\Delta_{x,y,z}v=f-\Delta_{x,y,z}B_0
\quad\text{in }\Omega,
\qquad v|_{\partial\Omega}=0.
\]
Linear Dirichlet theory gives a unique smooth solution with smooth
parameter dependence. For $m\ge2$, Corollary~\ref{cor:three} gives
the unique smooth admissible solution $u$.
In the comparisons below, we use Lemma~\ref{lem:comparison} for
$m\ge2$ and the maximum principle for $\Delta_{x,y,z}$ when $m=1$.

Since
\[
\det H_3(B_0+br)\ge f=\det H_3(u)\quad\text{in }\Omega,
\qquad B_0+br=u=B_0\quad\text{on }\partial\Omega,
\]
we have, $B_0+br\le u$, and then $br\le v$.
Set
\[
\delta=\min_{\overline\Omega}\bigl(f-\det H_3(B_0)\bigr)>0.
\]
By continuity, choose $0<b'\le b$ sufficiently small that
\[
\bigl\|\det H_3(B_0+b'r)-\det H_3(B_0)\bigr\|_{C^0(\overline\Omega)}
<\frac\delta2.
\]
Then $H_3(B_0+b'r)>0$ and
\[
\det H_3(B_0+b'r)<f-\frac\delta2<\det H_3(u),
\qquad B_0+b'r=u\quad\text{on }\partial\Omega.
\]
By comparison, $u\le B_0+b'r$, hence $v\le b'r$.
The choices of $b,b'$ are locally uniform in smooth families with a positive
gap, and smooth dependence in the nonlinear case follows from
Remark~\ref{prop:families}.

\begin{remark}
The same argument applies to any smooth $\Omega\Subset Y$ admitting a smooth
defining function $r$ with $H_3(r)>0$ on $\overline\Omega$, under the same
assumptions on $B_0$ and $f$.

Without the strict gap, the upper bound in \eqref{eq:two-sided} may fail:
if $H_3(B_0)>0$ and $f=\det H_3(B_0)$, uniqueness in
Corollary~\ref{cor:three} gives $v=0$, whereas $b'r<0$ in $\Omega$.
\end{remark}
\section{Periodic annuli}\label{sec:examples}

Let $0<a<b$ and $\TT^m=\RR^m/\ZZ^m$. Set
\[
\Omega_0=\{(x,y,z)\in\TT^m\times\RR^{2m}:
                     a^2<|y|^2+|z|^2<b^2\},
\qquad u_0=|y|^2+|z|^2.
\]
Formula \eqref{eq:H3} gives $H_3(u_0)=4I_m$.
Thus, for every smooth $0<f\le4^m$, Corollary~\ref{cor:three} gives a
unique smooth admissible solution with boundary values $a^2$ and $b^2$
on the inner and outer components, respectively.

The annuli $\Omega_0$ do not admit strictly
plurisubharmonic defining functions.
Indeed, lift to $\Omega_0\times\TT^m$ and use
$\rho=a^2-u_0$ near the inner boundary. Then $H(\rho)=-4I_m$.
Any other defining function with the same interior sign has the form
$r=c\rho$, where $c>0$ near the boundary. For a quaternionic tangent
column $\xi$, the product rule gives
\[
\xi^*H(r)\xi=c\xi^*H(\rho)\xi=-4c|\xi|^2
\quad\text{on the inner boundary}.
\]
Taking $\xi\ne0$ shows that neither a strictly
plurisubharmonic defining function on the lift nor a defining function with
$H_3(r)>0$ on $\overline\Omega_0$ can exist.

For a constant density $f=4^m c$, $0<c<1$, the solution can be written
explicitly. Put $s=|y|^2+|z|^2$ and seek $u=U(s)$. The chain rule gives
\[
\frac{\partial^2 U(s)}{\partial\bar q_a\,\partial q_b}
=4U'(s)\delta_{ab}+U''(s)s_{\bar a}s_b,
\qquad \sum_{a=1}^m|s_{\bar a}|^2=4s.
\]
The resulting matrix has eigenvalues $4U'$ with multiplicity $m-1$ and
$4(U'+sU'')$. Thus the equation becomes
\[
(U')^{m-1}(U'+sU'')=c.
\]
Integrating yields
\begin{equation}\label{eq:radial}
U(s)=a^2+\int_{a^2}^s(c+At^{-m})^{1/m}\,dt,
\end{equation}
where $A>0$ is chosen to satisfy
\[
\int_{a^2}^{b^2}(c+At^{-m})^{1/m}\,dt=b^2-a^2.
\]
The integral is strictly increasing in $A$, is smaller than $b^2-a^2$ at
$A=0$, and tends to infinity. Hence $A$ exists and is unique.
Finally, $U'>0$ and $U'+sU''=c/(U')^{m-1}>0$, so \eqref{eq:radial}
is admissible.

%
%
%

\section{Densities depending on the solution}\label{sec:solution-density}

\begin{theorem}\label{thm:solution-density}
Let $M\Subset X$ be as in Theorem~\ref{thm:main}, and let
$f\in C^\infty(\overline M\times\RR)$ satisfy $f>0$ and
$f_s=\partial f/\partial s\ge0$.
Suppose that $\underline u\in C^\infty(\overline M)$ satisfies
\[
H(\underline u)>0,\qquad
\det H(\underline u)\ge f(q,\underline u),
\qquad \underline u|_{\partial M}=\phi.
\]
Then there is a unique $u\in C^\infty(\overline M)$ satisfying
\begin{equation}\label{eq:solution-problem}
\det H(u)=f(q,u),\qquad H(u)>0\quad\text{in }M,
\qquad u|_{\partial M}=\phi.
\end{equation}
\end{theorem}

\begin{proof}
We replace \eqref{eq:path} by
\begin{equation}\label{eq:solution-path}
\begin{aligned}
\det H(u)&=F(t,q,u),\qquad u|_{\partial M}=\phi,\\
F(t,q,s)&=(1-t)\det H(\underline u)(q)+tf(q,s),
\qquad 0\le t\le1.
\end{aligned}
\end{equation}
Then $F>0$, $F_s\ge0$, and
$\det H(\underline u)\ge F(t,q,\underline u)$.
The path starts at $u=\underline u$.

\smallskip
\noindent\emph{Comparison and first derivatives.}
Let $L$ be given by \eqref{eq:linearization}, evaluated at $u$.
{By Lemma~\ref{lem:moore-variations}, we have}
\begin{equation}\label{eq:solution-comparison}
\begin{aligned}
L(\underline u-u)
&\ge\log\det H(\underline u)-\log\det H(u)\\
&\ge\log F(t,q,\underline u)-\log F(t,q,u)
 =c(q)(\underline u-u),
\end{aligned}
\end{equation}
where
\[
c(q)=\int_0^1\frac{F_s}{F}
 \bigl(t,q,u+\theta(\underline u-u)\bigr)\,d\theta\ge0.
\]
For each fixed smooth admissible $u$, $L$ is uniformly elliptic on
$\overline M$. Since $-c\le0$, the maximum principle yields
\[
(L-c)(\underline u-u)\ge0\ \text{in }M,\qquad
(\underline u-u)|_{\partial M}=0
\quad\Longrightarrow\quad \underline u\le u.
\]
With $h$ as in Section~\ref{sec:estimates},
\[
\Delta(u-h)=\tr H(u)>0\ \text{in }M,\qquad
(u-h)|_{\partial M}=0
\quad\Longrightarrow\quad u\le h.
\]
The boundary equality then bounds the normal derivative as in
\eqref{eq:boundary-C1}. Thus
\begin{equation}\label{eq:solution-C0}
\underline u\le u\le h,
\qquad
\|u\|_{C^0(\overline M)}+\sup_{\partial M}|\nabla u|\le C.
\end{equation}
Set $I=[\min_{\overline M}\underline u,\max_{\overline M}h]$ and
\[
F_-:=\min_{[0,1]\times\overline M\times I}F>0,
\qquad F_+:=\max_{[0,1]\times\overline M\times I}F.
\]
All partial derivatives of $\log F$ up to order two are uniformly bounded
on $[0,1]\times\overline M\times I$. Constants below
depend only on the fixed data, not on $u$ or $t$.

The full linearization of the equation is
\begin{equation}\label{eq:solution-linearization}
\mathcal L v=Lv-(\partial_s\log F)v,
\qquad \partial_s\log F\ge0.
\end{equation}
Here and below, the derivatives of $F$ are evaluated at $(t,q,u(q))$.
For a parallel real unit derivative $D$ as in
Lemma~\ref{lem:first-derivatives}, $D_q$ denotes differentiation in $q$
with $t$ and $s$ fixed. Differentiating the equation gives
\[
L(Du)=D_q\log F+(\partial_s\log F)Du,
\qquad \mathcal L(Du)=D_q\log F.
\]
Choose $C_0=1+\max_{\overline M}\underline u$.
By \eqref{eq:basic-L},
\begin{equation}\label{eq:solution-test}
\begin{aligned}
\mathcal L(\underline u-C_0)
&=L\underline u+(\partial_s\log F)(C_0-\underline u)\\
&\ge\kappa\sum_{a=1}^m u^{a\bar a}
 \ge\kappa mF_+^{-1/m}>0.
\end{aligned}
\end{equation}
Hence, for a fixed sufficiently large $A$,
\[
\mathcal L\bigl(\pm Du+A(\underline u-C_0)\bigr)
\ge-C+A\kappa mF_+^{-1/m}\ge0.
\]
The maximum principle and \eqref{eq:solution-C0} yield
\begin{equation}\label{eq:solution-C1}
\sup_M|Du|\le\sup_{\partial M}|Du|
             +A\|C_0-\underline u\|_{C^0},
\qquad \|u\|_{C^1(\overline M)}\le C.
\end{equation}

\smallskip
\noindent\emph{Reduction of second derivatives to the boundary.}
The calculation in \eqref{eq:second-differentiated} gives
$L(D^2u)\ge D^2[\log F(t,q,u(q))]$. Therefore
\begin{equation}\label{eq:solution-second-differentiation}
\begin{aligned}
\mathcal L(D^2u)
&\ge D^2[\log F(t,q,u(q))]-(\partial_s\log F)D^2u\\
&=D_q^2\log F+2(D_q\partial_s\log F)Du
                  +(\partial_s^2\log F)(Du)^2
\ge-C,
\end{aligned}
\end{equation}
where the last inequality uses \eqref{eq:solution-C1}.
Increasing $A$ if necessary, \eqref{eq:solution-test} gives
\[
\mathcal L\bigl(D^2u+A(\underline u-C_0)\bigr)
\ge-C+A\kappa mF_+^{-1/m}\ge0.
\]
Since the zeroth-order coefficient of $\mathcal L$ is nonpositive,
the maximum principle implies
\begin{equation}\label{eq:solution-C2-reduction}
\sup_M D^2u
\le\max\left\{0,\sup_{\partial M}D^2u\right\}
   +A\|C_0-\underline u\|_{C^0}.
\end{equation}

\smallskip
\noindent\emph{Boundary estimates and existence.}
Use the local inverses of $\pi:\HH^m\to X$ from
Section~\ref{sec:normal}. On each lifted boundary neighborhood,
\eqref{eq:solution-path} has the same form, with $F>0$, $F_s\ge0$,
and the same subsolution.
For $v=u-\underline u$ and any real coordinate $r$,
\[
Lv_r=(\log F)_{q_r}+(\partial_s\log F)u_r-L\underline u_r,
\qquad |Lv_r|\le C\sum_a u^{a\bar a},
\]
by \eqref{eq:solution-C1} and \eqref{eq:basic-L}. Hence the barrier
and derivative-square estimates of Section~\ref{sec:normal} remain valid.
The local boundary estimates of Ma-Zhu
\cite[Section 3, Step 3, and Section 4]{Zhu} apply directly
to \eqref{eq:solution-path}. The constants in the boundary estimates
are controlled by
\eqref{eq:solution-C1}, the data on $\overline M\times I$, and the
fixed boundary charts. The artificial spherical boundary is controlled
by the same $C^1$ bound. Thus
{
\begin{equation}\label{eq:solution-boundary-C2}
\sup_{p\in\partial M}\sup_{|D|=1}|D^2u(p)|\le C.
\end{equation}
\par}
Combining \eqref{eq:solution-boundary-C2} with
\eqref{eq:solution-C2-reduction} and $\Delta u>0$ gives
{
\begin{equation}\label{eq:solution-global-C2}
\|u\|_{C^2(\overline M)}\le C,
\qquad \lambda_-I_m\le H(u)\le\lambda_+I_m,
\end{equation}
\par}
as at the end of the proof of Proposition~\ref{prop:C2}.
Using \eqref{eq:solution-global-C2}, the higher-order estimates
and closedness follow as in Section~\ref{sec:existence}.
Linear Dirichlet theory gives the isomorphism
\[
\mathcal L:C^{2,\alpha}_0(\overline M)
             \longrightarrow C^\alpha(\overline M),
\]
so the implicit function theorem gives openness.
Openness and closedness yield existence at $t=1$.
Applying \eqref{eq:solution-comparison}
to two solutions in both orders gives uniqueness.

\end{proof}


\begin{thebibliography}{99}

\bibitem{Alesker1}
S.~Alesker,
\emph{Non-commutative linear algebra and plurisubharmonic functions of
quaternionic variables},
Bull. Sci. Math. \textbf{127} (2003), 1--35.

\bibitem{Alesker2}
S.~Alesker,
\emph{Quaternionic Monge--Amp\`ere equations},
J. Geom. Anal. \textbf{13} (2003), 205--238.

\bibitem{Alesker2013}

S.~Alesker,
\emph{Solvability of the quaternionic Monge--Amp\`ere equation on
compact manifolds with a flat hyperK\"ahler metric},
Adv. Math. \textbf{241} (2013), 192--219.

\bibitem{AleskerShelukhin2013}

S.~Alesker and E.~Shelukhin,
\emph{On a uniform estimate for the quaternionic Calabi problem},
Israel J. Math. \textbf{197} (2013), 309--327.

\bibitem{AleskerVerbitsky2010}

S.~Alesker and M.~Verbitsky,
\emph{Quaternionic Monge--Amp\`ere equation and Calabi problem
for HKT-manifolds},
Israel J. Math. \textbf{176} (2010), 109--138.


\bibitem{CNS1}

L.~Caffarelli, L.~Nirenberg and J.~Spruck,
\emph{The Dirichlet problem for nonlinear second-order elliptic
equations. I. Monge--Amp\`ere equation},
Comm. Pure Appl. Math. \textbf{37} (1984), 369--402.

\bibitem{CKNS}

L.~Caffarelli, J.~J. Kohn, L.~Nirenberg and J.~Spruck,
\emph{The Dirichlet problem for nonlinear second-order elliptic
equations. II. Complex Monge--Amp\`ere, and uniformly elliptic, equations},
Comm. Pure Appl. Math. \textbf{38} (1985), 209--252.




\bibitem{ChengYau1976}
S.-Y.~Cheng and S.-T.~Yau,
\emph{On the regularity of the solution of the $n$-dimensional
Minkowski problem},
Comm. Pure Appl. Math. \textbf{29} (1976), 495--516.





\bibitem{Evans}
L.~C. Evans,
\emph{Classical solutions of fully nonlinear, convex, second-order
elliptic equations},
Comm. Pure Appl. Math. \textbf{35} (1982), 333--363.

\bibitem{GMZComplexExterior}

Z.~Gao, X.-N.~Ma and D.~Zhang,
\emph{The exterior Dirichlet problem for the homogeneous complex
$k$-Hessian equation},
Adv. Nonlinear Stud. \textbf{23} (2023), no.~1, Article 20220039, 33~pp.


\bibitem{GMZRealPunctured}

Z.~Gao, X.-N.~Ma and D.~Zhang,
\emph{The Dirichlet problem of the homogeneous $k$-Hessian equation
in a punctured domain},
preprint (2023),
\href{https://arxiv.org/abs/2303.07976}{arXiv:2303.07976}.

\bibitem{GMZComplexPunctured}

Z.~Gao, X.-N.~Ma and D.~Zhang,
\emph{The Dirichlet problem of homogeneous complex $k$-Hessian equation
in a $(k-1)$-pseudoconvex domain with isolated singularity},
preprint (2023),
\href{https://arxiv.org/abs/2304.08407}{arXiv:2304.08407}.

\bibitem{GMYZ2026}

Z.~Gao, X.-N.~Ma, B.~Yu and D.~Zhang,
\emph{$C^{1,1}$ regularity global estimates for homogeneous complex
Hessian equations on punctured domains},
preprint (2026),
\href{https://arxiv.org/abs/2609.14692}{arXiv:2609.14692}.

\bibitem{GentiliVezzoni2026}

G.~Gentili and L.~Vezzoni,
\emph{Fully non-linear elliptic equations on compact hyperk\"ahler
manifolds},
Calc. Var. Partial Differential Equations \textbf{65} (2026), Article~91.


\bibitem{GT}
D.~Gilbarg and N.~S. Trudinger,
\emph{Elliptic Partial Differential Equations of Second Order},
Classics in Mathematics, Springer, 2001.

%


\bibitem{Guan}
B.~Guan,
\emph{The Dirichlet problem for complex Monge--Amp\`ere equations and
regularity of the pluri-complex Green function},
Comm. Anal. Geom. \textbf{6} (1998), 687--703.

\bibitem{Guan2009}
B.~Guan,
\emph{The Dirichlet problem for complex Monge--Amp\`ere equations and
applications},
in \emph{Trends in Partial Differential Equations},
Adv. Lect. Math. (ALM), vol.~10,
Higher Education Press, Beijing, and International Press, Somerville, MA,
2009, pp.~53--97.



\bibitem{GuanSpruck1993}
B.~Guan and J.~Spruck,
\emph{Boundary-value problems on $\mathbb S^n$ for surfaces of constant
Gauss curvature},
Ann. of Math. (2) \textbf{138} (1993), 601--624.





\bibitem{Krylov1}
N.~V. Krylov,
\emph{Boundedly nonhomogeneous elliptic and parabolic equations},
Math. USSR-Izv. \textbf{20} (1983), 459--492.

\bibitem{Krylov2}
N.~V. Krylov,
\emph{Boundedly nonhomogeneous elliptic and parabolic equations in a domain},
Math. USSR-Izv. \textbf{22} (1984), 67--97.





\bibitem{MaZhangExterior}

X.-N.~Ma and D.~Zhang,
\emph{The exterior Dirichlet problem for the homogeneous
$k$-Hessian equation},
preprint (2022), revised 2024,
\href{https://arxiv.org/abs/2207.13504}{arXiv:2207.13504}.

\bibitem{OpenAI}
OpenAI,
\emph{Local smoothing for the Euclidean wave equation},
preprint, October 25, 2026.

\bibitem{Po78}
A.V. Pogorelov, {\em The Minkowski multidimensional problem},
Wiley, New York, (1978).





\bibitem{Yau1978}
S.-T.~Yau,
\emph{On the Ricci curvature of a compact K\"ahler manifold and the
complex Monge--Amp\`ere equation, I},
Comm. Pure Appl. Math. \textbf{31} (1978), 339--411.


\bibitem{Zhu}
J.~Zhu,
\emph{Dirichlet problem of quaternionic Monge--Amp\`ere equations},
Israel J. Math. \textbf{214} (2016), 597--619.


\end{thebibliography}
\end{document}